\documentclass[11pt,leqno]{article}
\usepackage{amsmath,amssymb,amsfonts}
\usepackage{enumerate}
\usepackage[utf8]{inputenc}
\allowdisplaybreaks

\author{Mikhail {\sc Kamenskii}\footnote{
Department of Mathematics,
Voronezh State University, Voronezh, Russia.
E-Mail: Mikhail@kam.vsu.ru}, 
Omar {\sc Mellah}
\footnote{Normandie Univ., Laboratoire Rapha\"el Salem,
UMR CNRS 6085, Rouen, France 
and Department of Mathematics
Faculty of Sciences
University Mouloud Mammeri of Tizi-Ouzou, Algeria. 
E-Mail: omellah@yahoo.fr} 
and 
Paul {\sc Raynaud de Fitte}
\footnote{Normandie Univ., Laboratoire Rapha\"el Salem,
UMR CNRS 6085, Rouen, France.
E-Mail: prf@univ-rouen.fr
} 
}
\title{Weak Averaging of Semilinear Stochastic Differential Equations
  with Almost Periodic Coefficients}

\newtheorem{theo}{Theorem}[section]
\newtheorem{prop}[theo]{Proposition}

\newtheorem{lem}[theo]{Lemma}
\newtheorem{definition}[theo]{Definition}

\newtheorem{remark}[theo]{Remark}
\newenvironment{rem}{\begin{remark}\em}{\end{remark}}

\def\cprime{$'$}

\renewcommand\epsilon{\varepsilon}

\newcommand\R{\mathbb{R}}

\newcommand\un[1]{\,\rlap{{1}}\kern.22em \mbox{l}_{#1}} 
\newcommand\range[1]{\mathop{\text{\rm range}} (#1) }

\newcommand\proof{\noindent {\bf Proof}\ }
\newcommand\proofof[1]{\noindent{\bf Proof of {#1}}\ }
\renewcommand\square{\fbox{\rule{0em}{.3em}\rule{.3em}{0em}} \qquad}
\newcommand\finpr{\hfill$\square\qquad$\medskip\par}

\newcommand\tq{;\,} 

\newcommand\CCO[1]{\left( #1 \right)}
\newcommand\norm[1]{\left\Vert #1 \right\Vert}

\newcommand\accol[1]{\left\{#1\right\}}
\newcommand\scal[1]{\left\langle #1\right\rangle}

\newcommand\expect{\mathop{\text{\rm E}}\nolimits}

\newcommand\loi{\mathop{\text{\rm law}}}
\newcommand\law[1]{ \loi({#1}) } 

\newcommand\ellp[1]{\mathop{\text{\rm L}}\nolimits^{#1}}

\newcommand\probas[1]{{\mathcal{P}}\CCO{#1}}

\newcommand\trace{\mathop{\mbox{\rm tr}}}
\newcommand\Dom{\mathop{\mbox{\rm Dom}}}

\newcommand\cobar{\mathop{\overline{\text{\rm co}}}} 

\newcommand\esp{{\mathbb E}} 
\newcommand\distE{d_\esp}

\newcommand\polish{\mathbb{U}}

\newcommand\esprob{\Omega}
\newcommand\tribu{\mathcal{F}}
\newcommand\prob{\mathop{\text{\rm P}}\nolimits}

\newcommand\transl[1]{\widetilde{#1}}

\newcommand\C{\mathcal{C}}
\newcommand\FF{\mathcal{K}}
\newcommand\h{\mathbb{H}}
\newcommand\NN{\mathcal{N}}
\newcommand\V{\mathbb{V}}

\newcommand\CUB{\mbox{\rm CUB}}
\newcommand\BL{\mbox{\rm BL}}
\newcommand\bl{\mbox{\tiny\rm BL}}
\newcommand\WASS{\mathop{\mbox{\rm\bf W}}\nolimits^2}
\newcommand\wass{\mathop{\mbox{\rm\bf W}^2_{[\Tau,T]}}}

\newcommand\BDG{C_p} 
\newcommand\bdg[1]{C_{#1}}
\newcommand\convol{\mathfrak{C}_{p,T-\Tau}}
\newcommand\Tau{\kappa}
\newcommand\sig{s}

\newcommand\GG{\mathfrak{G}}

\newcommand\Bspace{\mathbb{E}}

\begin{document}

\maketitle

\begin{abstract}
An averaging result is proved for stochastic evolution equations
with highly oscillating coefficients. This result applies in
particular to equations with almost periodic coefficients. 
The convergence to the solution to the averaged equation 
is obtained in distribution, 
as in previous works by Khasminskii and Vrko\v{c}.

This version corrects two minor errors from our paper published in 
J. Math. Anal. Appl. 427(1):336–364, 2015.
\end{abstract}


\section{Introduction}
Since the classical work of N.M.~Krylov and N.N.~Bogolyubov
\cite{krylov-bogoliubov43book} devoted to the analysis, by  
the method of averaging, of the problem  
of the dependence on a small parameter $\varepsilon > 0$ of almost periodic
solutions to ordinary differential equation containing  
terms of frequency of order $\frac{1}{\varepsilon}$, several articles
and books have appeared, which develop this method for different kinds of 
differential equations. 
See the bibliography in the
book of V.Sh.~Burd \cite{BURD2007}, 
where a list of books related to this problem for deterministic
differential equations is presented. 
We note here that
 the authors of these papers 
are greatly influenced by the books of N.N.~Bogolyubov and 
A.Yu.~Mitropolskii \cite{bogolyubov-mitropolsky61book}
and 
M.A.~Krasnosel{\cprime}ski{\u\i}, V.Sh.~Burd and Yu.S.~Kolesov
\cite{Kraskoselski-Burd-Kolesov73book}. 

The method of averaging has been applied of course to stochastic
differential equations, but in general it was applied  
 to the initial problem in a finite interval, see for
 example \cite{Khasminskii1968}. Even in this case we  
can see a great 
difference with the deterministic case. 
To ensure the strong convergence in a space of stochastic processes, we must 
assume such convergence of the stochastic
term when $\varepsilon\rightarrow 0$, 
which virtually excludes the 
consideration of 
high frequency oscillation of this term. 
R.Z.~Khasminskii \cite{Khasminskii1968} 
has shown, in a finite dimensional setting, 
that it is possible to overcome this 
problem if one only looks for convergence {\em in distribution} to the
solution to the averaged equation. 
Later Ivo Vrko\v{c} \cite{vrkoc95weakaveraging} generalized this
result in a Hilbert space setting, 
for which the initial problem was at this time already well
developed 
(see for example the book of Da Prato and Zabczyck \cite{dapratozabczyk92book}).

During the last 20 years an intensive study of the problem of
existence of almost periodic solutions to stochastic  
differential equations was performed by L.~Arnold, C.~Tudor, G.~Da
Prato   
(see in particular \cite{DaPrato-Tudor95,Arnold-Tudor98})
and later by P.H.~Bezandry and T.~Diagana 
\cite{bezandry-diagana07,bezandry-diagana07quadratic,%
bezandry-diagana09quadratic}. 
For the first group, an almost periodic solution 
 means that the stochastic process generates an almost periodic
 measure on the paths space. 
The second group claims the 
existence of square mean almost periodic solutions, but square mean
almost periodicity seems to be a too strong  
property for solutions to SDEs, see counterexamples in
\cite{mellah-prf2012contrexemples}. 
 
In this paper we propose the averaging principle for 
solutions  
 to a family of semilinear stochastic differential equations in Hilbert
 space
which are almost periodic in distribution. 
The second member of these equations contains a 
 high frequency term. Under the Bezandry-Diagana conditions, we
 establish the convergence in distribution (actually, in Wasserstein distance)
of the solutions to  
these equations to the solution to the averaged equation in the
sense of Khasminskii-Vrko\v{c}, 
with a weaker hypothesis than \cite{vrkoc95weakaveraging} on the
linear evolution semigroup. 

The paper is organized as follows: 
The next section is devoted to the
notations and preliminaries. 
We then prove in Section \ref{sect:Form.prob} 
that the solutions to the equations we consider are
 almost periodic in distribution, 
when their coefficients are almost periodic.  
In section \ref{sect:averaging}, 
we prove the 
fondamental averaging result of this paper. 

This version corrects two minor errors from the paper published in 
J. Math. Anal. Appl. 427(1):336–364, 2015: 
on the one hand, 
the value of the coefficient $\theta$ in Theorem \ref{theo:main} is
slightly different, on the other hand the argument given for the proof
of the convergence of $I_2$ and $I_4$ in the second step of the proof
of Theorem \ref{theo:main} has been corrected.

\section{Notations and Preliminaries}
In the sequel, $(\h_1, \|.\|_{\h_1})$ and $(\h_2, \|.\|_{\h_2})$ denote separable
Hilbert spaces and $L(\h_1, \h_2)$ (or $L(\h_1)$ if 
 $\h_1= \h_2$) is the space of all bounded linear operators from
 $\h_1$ to $\h_2$, whose norm will be denoted by  
$\|.\|_{L(\h_1, \h_2)}$. 
If $A \in L(\h_1)$ then $A^*$ denotes its
adjoint operator and if $A$ is a nuclear operator, 
$$|A|_\NN = \sup\accol{\sum_i|<Ae_i, f_i>|, \{e_i\}, \{f_i\}
  \text{ orthonormal bases of }  \h_1}$$ is  
the nuclear norm of $A$.

\subsection{Almost periodic functions}
Let $(\esp, d)$ be a separable metric space, we denote by $C_b(\esp)$ the
Banach space of continuous and bounded functions  
$f:\esp\rightarrow\R$ with $\parallel f \parallel_\infty = \sup_{x\in
  \esp}|f(x)|$ and by $\probas{\esp}$ the set of  
all probability measures onto $\sigma$-Borel field of $\esp$. For
$f\in C_b(\esp)$ we define 
\begin{gather*}
\parallel f \parallel_L 
= \sup\Bigl\{\dfrac{f(x)-f(y)}{d_\esp(x,y)}: x \neq y\Bigl\}\\ 
\parallel f \parallel_{\bl} = \max\{\parallel
f \parallel_\infty,\parallel f \parallel_L \}
\end{gather*}
and we define  
$$\BL(\esp) 
= \bigl\{f\in C_b(\esp);\parallel f \parallel_{\bl}<\infty \bigl\}.$$
For $\mu, \nu \in \probas{\esp}$ we define 
$$d_{\bl}(\mu, \nu) = \sup_{\parallel f \parallel_{\bl}\leq 1}\Bigl|
\int_\esp fd(\mu-\nu)\Bigl|  $$  
which is a complete metric on $\probas{\esp}$ and generates the narrow
(or weak) topology, i.e.~the coarsest topology on $\probas{\esp}$ such
that the mappings $\mu\mapsto\mu(f)$ are continuous for all bounded
continuous $f :\,\esp\rightarrow\R$.

Let $(\esp_1, d_1)$ and $(\esp_2, d_2)$ be separable and complete
metric spaces. Let $f$ be a continuous 
mapping from $\R$ to $\esp_2$ 
(resp. from $\R \times \esp_1 $ to $\esp_2$). 
Let $\mathcal{K}$ be a set of subsets of $\esp_1$. 
The function $f$ is
said to be {\em almost periodic} 
(respectively {\em almost periodic uniformly 
with respect to $x$ in elements of $\mathcal{K}$})
if for every  $\varepsilon>0$ 
(respectively for every  $\varepsilon>0$
and every subset $K\in\mathcal{K}$), 
there exists a constant $l(\varepsilon, K)>0$ 
such that any interval of length $l(\varepsilon, K)$ contains at least
a number $\tau$ for which  
\begin{gather*}
\sup_{t\in\R} d_2( f(t+\tau) , f(t))
<\varepsilon\\
\text{ (respectively }\sup_{t\in\R} \sup_{x\in K} d_2( f(t+\tau, x) , f(t, x))
<\varepsilon\text{).} 
\end{gather*}

A characterization of almost periodicity is given in the following
result, due to Bochner:
\begin{theo}\label{theo:characterisation}
 Let $f: \R \rightarrow \h_1 $ be continuous. Then the
 following statements are equivalent 
\begin{itemize}
 \item $f$ is almost periodic.
\item The set of translated functions $\{ f(t+.)\}_{t \in \R}$ is relatively
  compact in $\C(\R; \esp_2)$ with respect to the uniform norm.
\item $f$ satisfies {\em Bochner's double sequence criterion}, 
that is, for every pair of sequences 
$\{\alpha^{'}_n\} \subset \R$ and $\{ \beta^{'}_n\} \subset \R$, there are subsequences 
$(\alpha_n)\subset (\alpha'_n)$ and
  $(\beta_n)\subset (\beta'_n)$ respectively with same indexes such
  that, for every $t\in \R$, the limits  
\begin{equation}\label{def:Bochner double sequence}
 \lim_{n\rightarrow\infty}\lim_{m\rightarrow\infty}f(t+\alpha_n+\beta_m)
 \text{ and }\lim_{n\rightarrow\infty}f(t+\alpha_n+\beta_n), 
\end{equation}
exist and are equal.
\end{itemize}
 
\end{theo}
\begin{rem} \label{Boch:advantage}
\rule{1em}{0em}
\begin{enumerate}[{\rm (i)}]
\item A striking property of Bochner's double sequence criterion is
   that the limits in (\ref{def:Bochner double sequence}) exist in any
   of the three modes of convergences: pointwise, uniform on compact
   intervals and uniform on $\R$ (with respect to $\distE$). This
   criterion has thus the avantage that it allows to establish uniform
   convergence by checking pointwise convergence. 
\item The previous result holds for the metric spaces  $(\probas{\esp}, d_{\bl})$\\
 and 
 $(\probas{\C (\R, \esp)}, d_{\bl})$  
 
\end{enumerate}
\end{rem}

\subsection{Almost periodic stochastic processes}
Let $(\esprob, \tribu, \prob)$ be a probability space. Let $X:\R
\times \esprob \rightarrow \h_2$ be a stochastic  
process. 
We denote by $\law{X(t)}$
the distribution of the random variable $X(t)$.  
Following Tudor’s terminology \cite{Tudor95ap_processes}, we say that
X has {\em almost periodic one-dimensional distributions} 
if the mapping $t\mapsto \law{X(t)}$ from
$\R$ to $(\probas{\h_2}, d_{\bl})$ is almost periodic. 

If $X$ has continuous trajectories, 
we say that $X$ is {\em almost periodic in distribution} 
if the mapping 
$t\mapsto\law{X(t+.)}$ 
from $\R$ to $\probas{C(\R;\h_2)}$ is almost periodic, 
where 
$C(\R;\h_2)$ is endowed 
with the uniform convergence on compact intervals 
and $\probas{C(\R;\h_2)}$ is endowed with the distance $d_{\bl}$.

Let $ \ellp{2}(\prob, \h_2)$ be the space of $\h_2$-valued 
random variables with a finite quadratic-mean. 
We say that a stochastic process $X: \R\rightarrow \ellp{2}(\prob, \h_2)$ is  
{\em square-mean continuous} if, for every $s\in\R$, 
$$\lim_{t\rightarrow s} \expect\|X(t) - X(s)\|^2 _ {\h_2} =0.$$ 
We denote by $\CUB\bigl(\R,\ellp{2}(\prob, \h_2)\bigl)$
 the Banach space of square-mean continuous and uniformly bounded stochastic
 processes, endowed with the norm  
$$\|X\|^2_{\infty}= \sup_{t\in\R}(\expect\|X(t)\|^2 _ {\h_2}).$$
A square-mean continuous stochastic process 
$X: \R\rightarrow \ellp{2}(\prob, \h_2)$ is
said to be {\em square-mean almost periodic} if, for each  
$\varepsilon>0$, there exists $l(\varepsilon) >0$ 
 such that any interval of length $l(\varepsilon)$ contains at least a
 number $\tau$ for which  
$$\sup_{t\in\R} \expect\|X(t+\tau) - X(t)\|^2 _ {\h_2}<\varepsilon.$$ 

The next theorem is interesting in itself, but we shall not use it in
the sequel. 
\begin{theo}\label{theo:Sq-AP}
Let $F: \R \times \h_2 \rightarrow \h_2$ be an almost periodic function
uniformly with respect to $x$ in compact subsets  
of $\h_2$ such that 
$$\|F(t,x)\|_{\h_2}\leq C_1(1 + \|x\|_{\h_2}) \text{ and }
\|F(t,x)-F(t,y)\|_{\h_2}\leq C_2\| x-y\|_{\h_2}. $$ 
Then the function 
$$\transl{F}: \R\times \ellp{2}(\prob, \h_2) \rightarrow \ellp{2}(\prob, \h_2)$$
(where $\transl{F}(t, Y)(\omega) = F(t, Y(\omega))
\text{ for every }\omega \in \esprob$) is square-mean almost periodic 
 uniformly with respect to $Y$ in compact subsets of $\ellp{2}(\prob, \h_2)$.
\end{theo}
\proof
Let us prove that for each $Y\in \ellp{2}(\prob, \h_2)$ the process 
$\transl{F}_Y: \R \rightarrow \ellp{2}(\prob, \h_2)$, 
 $t\mapsto\transl{F}(t, Y)$ is almost periodic.

For every $\delta>0$, there exists a compact subset $S$ of $\h_2$ such
that $$\prob\{Y\notin S\}\leq \delta.$$  
Let $\varepsilon>0$, then there exist $\delta>0$ and a compact subset
$S$ of $\h_2$ such that  
$\prob\{Y\notin S\}\leq \delta$ and $$\int_{\{Y\notin S\}}
\bigl(1+ \|Y\|^2_{\h_2}\bigl)d\prob< \frac{\varepsilon}{4C_1}.$$ Since
$F$ is almost periodic uniformly with respect to $x$ in the compact subset $S$, 
there exists a constant
$l(\varepsilon, S)>0$ such that any interval of  
length $l(\varepsilon, S)$ contains at least a number $\tau$ for
which $$\sup_t\|F(t+\tau, x) - F(t,x)\|_{\h_2}<  
\frac{\sqrt{ \varepsilon}}{\sqrt{2}} \, \, \text{ for all } x\in S.$$ 
We have 
\begin{align*}
\expect(\|F(t+\tau, Y) - F(t,Y)\|^2_{\h_2}) =
& \int_{\{Y\in S\}}\|F(t+\tau, Y) - F(t,Y)\|^2_{\h_2}d\prob \\
& + \int_{\{Y\notin S\}}\|F(t+\tau, Y) - F(t,Y)\|^2_{\h_2}d\prob\\
& < \frac{\varepsilon}{2} 
+ 2 C_1\int_{\{Y\notin S\}}\bigl(1+ \|Y\|^2_{\h_2}\bigl)d\prob\\
&  < \frac{\varepsilon}{2} + \frac{\varepsilon}{2} = \varepsilon. 
\end{align*}
Therefore the process $ \transl{F}_Y$ is almost periodic. 
Since $ \transl{F}$ is Lipschitz, it is almost periodic 
uniformly with respect to $Y$ in 
compact subsets of $\ellp{2}(\prob, \h_2)$ 
(see \cite[Theorem 2.10 page 25]{Fink74}).
\finpr
\begin{prop}\label{prop: mean}
Let $\mathcal{K}$ be a set of subsets of $\h_2$. 
Let $F: \R \times \h_2 \rightarrow \h_2$, $(t,x)\mapsto F(t,x)$, 
be almost periodic, uniformly
with respect to $x$ in elements of $\mathcal{K}$. 
There exists a continuous function $F_0: \h_2 \rightarrow \h_2$  
such that 
\begin{multline}\label{eq:moyenne_simple}
\lim_{t\rightarrow\infty}\frac{1}{t}\int_{\Tau}^{\Tau+t} F(s,x)\,ds=F_0(x)\\
\text{ for every }\Tau\in\R,
\text{ uniformly with respect to $x$ in elements of $\mathcal{K}$. } 
\end{multline}
Furthermore, if $F(t,x)$ is Lipschitz in $x\in \h_2$
 uniformly with respect to $t\in \R$, the mapping $F_0$ is
 Lipschitz too.
\end{prop}
\proof
Let $K\in\mathcal{K}$, and let $\mathbb{X}=\h_2^{K}$ be the Banach
space of all mappings from $K$ to $\h_2$, endowed with the supremum
norm. 
For every $t\in\R$, set 
$$\widehat{F}(t)=(F(t,x))_{x\in K},\ 
\widehat{F}_0=(F_0(x))_{x\in K}.$$
By \cite[Theorem 6.11]{Corduneanu68}, we have, for the norm of  $\mathbb{X}$,
$$\lim_{t\rightarrow\infty}\frac{1}{t}\int_{\Tau}^{\Tau+t}\widehat{F}(s)\,ds
   =\widehat{F}_0,$$
which proves  \eqref{eq:moyenne_simple}. 
Alternatively, on may use the proof of \cite[Theorem 3.1]{Fink74},
though the result is given in a finite dimensional setting and for a
compact set $K$. 

The Lipschitz property of $F_0$ is trivial. 
\finpr

Let $Q\in L(\h_1)$ be a linear operator. Then $Q$ is a bijection from 
$\range{Q}=Q(\h_1)$ to $(\ker Q)^\perp$. 
We denote by $Q^{-1}$ the {\em pseudo-inverse} of $Q$
(see \cite[Appendix C]{Prevot-Rockner07book} 
or \cite[Appendix B.2]{dapratozabczyk92book}), that is, the inverse of
the mapping $(\ker Q)^\perp\rightarrow\range{Q}$, $x\mapsto Q(x)$.  
Note that $\range{Q}$ is a Hilbert space for the scalar product 
$\scal{x,y}_{\range{Q}}=\langle Q^{-1}(x),Q^{-1}(y)\rangle$.

\begin{prop}\label{rem:mean}
Let $\mathcal{K}$ be a set of subsets of $\h_2$. 
 Let $G:\R\times\h_2\rightarrow L(\h_1, \h_2)$, 
$t\mapsto G(t,x)$,
 be almost periodic uniformly with respect to $x$ in elements of
 $\mathcal K$, 
 and let $Q\in L(\h_1)$ be a self-adjoint nonnegative 
operator. 
Let $\h_0=\range{Q^{1/2}}$, endowed with 
$\scal{ x,y }_{\range{Q^{1/2}}}
=\allowbreak\langle Q^{-1/2}(x),Q^{-1/2}(y)\rangle.$  
There exists a continuous function $G_0:\h_2\rightarrow L(\h_0, \h_2)$ such
that
\begin{multline}\label{eq:gaussassumption} 
  \lim_{t\rightarrow \infty}\bigl |\dfrac{1}{t}
  \int^{\Tau+t}_{\Tau}G(s,x)QG^*(s,x)\,ds - G_0(x)QG^{*}_{0}(x)
  \bigl|_{\NN} = 0\\ 
\text{ for all }\Tau\in\R,
\text{ uniformly with respect to $x$ in elements of $\mathcal{K}$, }  
 \end{multline}
where $G^*(s,x)=(G(s,x))^*$ and $G^{*}_{0}(x)=(G_{0}(x))^*$.
\end{prop}
\proof
Observe first that
$G_0(x)QG^{*}_{0}(x)=(G_0(x)Q^{1/2})(G_0(x)Q^{1/2})^*$, thus 
$G_0(x)$ does not need to be defined on the whole
space $\h_1$, it is sufficient that it be defined on $\h_0$.


Since $G$ is almost periodic, the function $H(s,x)=G(s,x)QG^*(s,x)$ is almost
periodic too, with positive self-adjoint nuclear values in $L(\h_2)$. 
Thus, reasoning as in the proof of Proposition \ref{prop: mean},
there exists a mapping 
$H_0 :\,\h_2\rightarrow L(\h_2)$
such that, for every $\Tau\in\R$, 
\begin{equation*}
\lim_{t\rightarrow\infty}\frac{1}{t}\int_{\Tau}^{\Tau+t} G(s,x)QG^*(s,x)\,ds=H_0(x)
\end{equation*}
uniformly with respect to $x$ in elements of $\mathcal{K}$.
By e.g.~\cite[Theorem 3.1]{Fink74}, $H_0$ is continuous. 
Thus
the mapping 
$$H_0^{1/2} : \,
\left\{\begin{array}{lcl}
\h_2&\rightarrow&L(\h_2)\\
x&\mapsto&(H_0(x))^{1/2}
\end{array}\right.$$
is continuous  
with positive self-adjoint values.

Let $G_0(x)=H_0^{1/2}(x)Q^{-1/2} :\,\h_0\rightarrow \h_2$.  
We then have, for every $x\in\h_2$, 
$$H_0(x)=G_0(x) Q (G_0(x))^*$$
and $G_0$ is continuous, 
which proves \eqref{eq:gaussassumption}.
\finpr
\section{Solutions almost periodic in distribution}%
\label{sect:Form.prob} 
We consider the semilinear stochastic differential equation, 
\begin{equation}\label{eq:SDE}
 dX_t = AX(t)dt + F(t, X(t))dt + G(t, X(t))dW(t), t\in\R
\end{equation}
Where $A: \Dom(A)\subset\h_2\rightarrow \h_2$ is a densely defined closed 
(possibly unbounded) linear operator, 
$F: \R\times\h_2 \rightarrow\h_2$, 
and $G: \R\times\h_2\rightarrow L(\h_1, \h_2)$ 
are continuous functions.
In this section, we assume that:
\newcounter{carbone}
\newcounter{carbon2}
\begin{enumerate}[{\rm (i)}]

 \item \label{cond:wienerpro}
 $W(t)$ is an $\h_1$-valued Wiener process with nuclear covariance operator 
 $Q$ (we denote by $\trace Q$ the trace of $Q$), 
 defined on a stochastic basis 
 $(\esprob,\tribu,(\tribu_t)_{t\in \R},\prob)$. 

\setcounter{carbon2}{\value{enumi}}  
 \item \label{cond:semigrp} $A : \Dom(A)\rightarrow \h_2$ is the
infinitesimal generator of  
a $C_0$-semigroup $(S(t))_{t\geq0}$ such that there exists a constant
$\delta>0$ with  
$$\|S(t)\|_{L(\h_2)}\leq e^{-\delta t}, t\geq0.$$
  \item \label{cond:croissance}
There exists a constant $K$ such that the mappings 
$F:\R\times\h_2\rightarrow \h_2$ and  
$G:\R\times\h_2\rightarrow L(\h_1, \h_2)$ 
satisfy
$$\|F(t,x)\|_{\h_{2}} + \|G(t,x)\|_{L(\h_1, \h_2)} \leq K (1+\|x\|_{\h_2}) $$
\item \label{cond:lipschitz}The functions $F$ and $G$ are Lipschitz,
  more precisely there exists a constant $K$ such that 
$$\|F(t,x) - F(t,y)\|_{\h_2} + \|G(t,x) - G(t,y)\|_{L(\h_1, \h_2)}\leq
K\|x - y\|_{\h_2}$$ for all $t\in \R$  
and $x, y \in \h_2$. 
\setcounter{carbone}{\value{enumi}}  
\item \label{cond:pp}
The mappings $F$ and  
$G$ are almost periodic in
$t\in\R$ 
uniformly with respect 
to $x$ in bounded subsets of $\h_2$. 

\end{enumerate}
The assumptions in the following theorem are contained in those of Bezandry and
Diagana 
\cite{bezandry-diagana07,bezandry-diagana07quadratic}. 
The result is similar to \cite[Theorem 4.3]{DaPrato-Tudor95}, with
different hypothesis and a different proof. 

\begin{theo}\label{theo:main}
 Let the assumptions (\ref{cond:wienerpro}) - (\ref{cond:pp}) be
 fulfiled and the constant  
$\theta=\dfrac{2K^2}{\delta}\CCO{\frac{1}{\delta}+\frac{\trace Q}{2}}<1$.
Then there exists a unique mild solution $X$ to (\ref{eq:SDE}) 
in $\CUB\bigl(\R,\ellp{2}(\prob, \h_2)\bigl)$. Furthermore, 
$X$ has a.e.~continuous
trajectories, and $X(t)$ can be explicitly
 expressed as follows, for each $t\in \R$:
\begin{equation}\label{eq:mildsol}
X(t) = \int^{t}_{-\infty}S(t-s)F\bigl(s, X(s)\bigr)ds + 
\int^{t}_{-\infty}S(t-s)G\bigl(s, X(s)\bigr)dW(s).
\end{equation}
If furthermore 
$\theta' =\dfrac{4K^2}{\delta}\CCO{\dfrac{1}{\delta} + \trace Q}<1$,
then $X$ is almost 
periodic in distribution. 
\end{theo}

To prove Theorem \ref{theo:main}, 
we need several preliminary results. 
Let us first recall the following result, which is given in a more general
form in \cite{DaPrato-Tudor95}:
\begin{prop}\label{prop:daprato-tudor3.1}%
(\cite[Proposition 3.1-(c)]{DaPrato-Tudor95})
Let $\tau\in\R$. 
Let $(\xi_n)_{0\leq n\leq\infty}$ 
be a sequence of square integrable $\h_2$-valued random variables. 
Let $(F_n)_{0\leq n\leq\infty}$ 
and $(G_n)_{0\leq n\leq\infty}$ be sequences of mappings from $\R\times\h_2$
to $\h_2$ and $L(\h_1, \h_2)$ respectively, satisfying 
\eqref{cond:croissance} and \eqref{cond:lipschitz} (replacing $F$ and
$G$ by $F_n$ and $G_n$ respectively, and the constant $K$ being
independent of $n$). 
For each $n$, let $X_n$ denote the solution to 
\begin{multline*}
X_n(t)=S(t-\tau)\xi_n\\
+\int^{t}_{\tau}S(t-s)F_n\bigl(s, X_n(s)\bigl)ds + 
\int^{t}_{\tau}S(t-s)G_n\bigl(s, X_n(s)\bigl)dW(s).
\end{multline*}
Assume that, for every $(t,x)\in \R\times\h_2$,  
\begin{gather*}
\lim_{n\rightarrow\infty}F_n(t,x)=F_{\infty}(t,x),\
\lim_{n\rightarrow\infty}G_n(t,x)=G_{\infty}(t,x),\\
\lim_{n\rightarrow\infty}d_{\bl}(\law{\xi_n,W},\law{\xi_\infty,W})=0, 
\end{gather*}
(the last equality takes place in $\probas{\h_2\times C(\R,\h_1)}$). 
Then we have in $\C([\tau,T]; \h_2)$, for any $T>\tau$,
$$\lim_{n\rightarrow\infty}d_{\bl}(\law{X_n},\law{X_\infty})=0.$$
\end{prop}

We need also a variant of Gronwall's lemma, taylored for mild
solutions. 
\begin{lem}\label{lem:gronwall}
Let $g :\R\to\R$ be a continuous function 
such that, 
for every $t\in\R$, 
\begin{equation}\label{hypothesegronwall}
0\leq g(t)\leq \alpha(t)
+\beta_1\int_{-\infty}^t e^{-\delta_1(t-s)}g(s)\,ds
+\dots
+\beta_n\int_{-\infty}^t e^{-\delta_n(t-s)}g(s)\,ds,
\end{equation}
for some locally integrable function $\alpha :\,\R\rightarrow\R$, 
and for some constants $\beta_1,\dots,\beta_n\geq 0$, 
and some constants $\delta_1,\dots,\delta_n>\beta$, 
where $\beta:=\sum_{i=1}^n\beta_i$. 
We assume that the integrals in the right hand
side of \eqref{hypothesegronwall} are convergent. 
Let $\delta=\min_{1\leq i\leq n}\delta_i$. 
Then, for every 
$\gamma\in ]0, \delta-\beta]$ 
such that $\int_{-\infty}^0 e^{\gamma s}\alpha(s)\,ds$ converges, we
have, for every $t\in\R$,
\begin{equation}\label{eq:gronwallinfinite}
g(t)\leq \alpha(t)+\beta\int_{-\infty}^t e^{-\gamma(t-s)}\alpha(s)\,ds. 
\end{equation}
In particular, if $\alpha$ is constant, we have 
\begin{equation}\label{eq:alphaconstant}
g(t)\leq \alpha\, \frac{\delta}{\delta-\beta}. 
\end{equation}
\end{lem}
\proof
Let $\beta'_i=\beta_i/\beta$,
$i=1,\dots,n$. 
We have 
\begin{multline*}
\dfrac{d}{dt}\CCO{e^{\gamma t}\sum_{i=1}^n
                    \beta'_i\int_{-\infty}^t e^{-\delta_i(t-s)}g(s)\,ds}\\
\begin{aligned}
=&\dfrac{d}{dt}\CCO{\sum_{i=1}^n
         \beta'_ie^{(\gamma-\delta_i)t}\int_{-\infty}^t e^{\delta_i s}g(s)\,ds}\\
=&\sum_{i=1}^n(\gamma-\delta_i)e^{(\gamma-\delta_i)t}
                     \beta'_i \int_{-\infty}^t e^{\delta_i s}g(s)\,ds
          +\sum_{i=1}^ne^{(\gamma-\delta_i) t}e^{\delta_i t}\beta'_i g(t)\\
=&e^{\gamma t} \CCO{
 g(t)+\sum_{i=1}^n(\gamma-\delta_i)
         \beta'_i\int_{-\infty}^t e^{-\delta_i(t-s)}g(s)\,ds
                      }\\
\leq& e^{\gamma t}\alpha(t).
\end{aligned}
\end{multline*}
The last inequality holds because $\gamma-\delta_i\leq -\beta$ and
$g\geq 0$. 
Integrating on $]-\infty,t]$, we get
\begin{align}
e^{\gamma t}\sum_{i=1}^n\beta'_i\int_{-\infty}^t e^{-\delta_i(t-s)} g(s)\,ds
&\leq \int_{-\infty}^t e^{\gamma s}\alpha(s)\,ds \notag\\
\intertext{(because both terms go to $0$ when $t\rightarrow-\infty)$, {\em i.e.}}
\sum_{i=1}^n\beta'_i\int_{-\infty}^t e^{-\delta_i(t-s)} g(s)\,ds
&\leq e^{-\gamma t} \int_{-\infty}^t e^{\gamma s}\alpha(s)\,ds. \label{etvoila} 
\end{align}
Using \eqref{etvoila} in \eqref{hypothesegronwall} yields
\begin{align*}
g(t)
&\leq \alpha(t)+\beta\sum_{i=1}^n\beta'_i\int_{-\infty}^t e^{-\delta_i(t-s)} g(s)\,ds
\leq \alpha(t)+\beta e^{-\gamma t}\int_{-\infty}^t e^{\gamma s}\alpha(s)\,ds.
\end{align*}
Inequality \eqref{eq:alphaconstant} is a direct consequence 
of \eqref{eq:gronwallinfinite}, 
with $\gamma=\delta-\beta$. 
\finpr

Lemma \ref{lem:gronwall} will help us (among other things) state a result 
(Proposition \ref{cor:Lpboundedness})
on estimation of  
the moments of solutions to \eqref{eq:mildsol}, which will be useful 
in the proof of Theorem \ref{theo:main} as well as
in the proof of Theorem \ref{theo:averaging}. 
To prove this result, 
we need also one of the moment inequalities for stochastic integrals due to
Novikov \cite{novikov1971moment_inequalities}. 
These inequalities are proved in a finite dimensional setting in
\cite{novikov1971moment_inequalities}, but their proofs extend easily
in infinite dimension. 
For the sake of completeness, we give the proof of the inequality we
need (for $p\geq 2)$, in our setting. 

\begin{lem}\label{lem:novikov}
Let $p\geq 2$, and let $Y$ be an $L(\h_1, \h_2)$-valued 
adapted stochastic process. 
We have, for every $t\geq 0$,  
\begin{equation}\label{eq:novikov}
\expect\norm{\int_0^t Y(s)\,dW(s)}^p
\leq \BDG 
\expect\CCO{ \int_0^t\trace\CCO{Y(s)QY^*(s)}\,ds }^{p/2}
\end{equation}
with
\begin{equation*}
\BDG=\frac{1}{(2c)^{p/2}}\CCO{\frac{2+2c}{p-1}-2^{p/2}}
\end{equation*}
for any $c>(p-1)2^{p/2-1}-1.$
In particular, $\bdg{2}=1$ (and in that case, \eqref{eq:novikov} is an
equality).
\end{lem}
\proof
We denote
$$Z(t)=\int_0^t Y(s)\,dW(s),\quad
V(t)=\int_0^t\trace\CCO{Y(s)QY^*(s)}\,ds.
$$

Let us first assume that 
$\norm{Y(t)}\leq M$ a.e.~for some constant $M$
and for every $t\geq 0$. 
Let $\alpha,c\geq 0$. 
Let 
 $$X(t)=\alpha+c V(t)+\norm{Z(t)}^2.$$
By Itô's formula (e.g.~\cite[Theorem 4.17]{dapratozabczyk92book}), 
we have
\begin{equation*}
dX(t)
=(1+c)\,dV(t)+2\scal{Z(t),Y(t)\,dW(t)}.
\end{equation*}
Denoting by $[X]$ the quadratic variation of $X$, we have thus
\begin{align*}
d\CCO{X(t)}^{p/2}
=&\frac{p}{2}X^{p/2-1}(t)\,dX(t)+\frac{1}{2}\frac{p}{2}\CCO{\frac{p}{2}-1}\,
X^{p/2-2}d[X](t)\\
=&\Bigg\lgroup\frac{p}{2}(1+c)\CCO{\alpha+c V(t)+\norm{Z(t)}^2}^{p/2-1}
          \trace\CCO{Y(t)QY^*(t)}\\
 &+\frac{p}{2}\CCO{p-2}
       \CCO{\alpha+c V(t)+\norm{Z(t)}^2}^{p/2-2}
         {\scal{Z(t),Y(t)Q^{1/2}}}^2
  \Biggr\rgroup\,dt\\
 &+p\CCO{\alpha+c V(t)+\norm{Z(t)}^2}^{p/2-1}
     \scal{Z(t),Y(t)\,dW(t)}.
\end{align*}
By the boundedness hypothesis on $Y$, the process 
$Z$ is a continuous martingale
and 
$e^{\norm{Z}-\frac{1}{2}V}$ is a supermartingale 
(actually it is a martingale, 
see \cite[Theorem IV.37.8]{rogers-williams2}).
We thus have 
\begin{equation*}
\expect\CCO{ e^{\norm{Z(t)}} }=\expect\CCO{ e^{\frac{1}{2}V(t)}
}\leq e^{\frac{1}{2}\trace{Q}M^2t},
\end{equation*}
thus the moments of any order of $\norm{Z(t)}$ are bounded. 
We deduce
\begin{multline}\label{eq:devlop}
\expect\CCO{{\alpha+c V(t)+\norm{Z(t)}^2}}^{p/2}\\
\begin{aligned}
=&\alpha^{p/2}
  +\frac{p}{2}(1+c)\expect\int_0^t \CCO{\alpha+c V(s)+\norm{Z(s)}^2}^{p/2-1}
          \trace\CCO{Y(s)QY^*(s)}\,ds\\
 &+\frac{p}{2}\CCO{p-2}\expect\int_0^t
      \CCO{\alpha+c V(s)+\norm{Z(s)}^2}^{p/2-2}
         \scal{Z(s),Y(s)Q^{1/2}}^2\,ds.   
\end{aligned}
\end{multline}
In particular, for $c=0$, 
\begin{multline*}
\expect\CCO{{\alpha+\norm{Z(t)}^2}}^{p/2}\\
\begin{aligned}
\leq &\alpha^{p/2}
      +\frac{p}{2}\expect\int_0^t \CCO{\alpha+\norm{Z(s)}^2}^{p/2-1}
          \trace\CCO{Y(s)QY^*(s)}\,ds\\
     &+\frac{p}{2}\CCO{p-2}\expect\int_0^t
      \CCO{\alpha+\norm{Z(s)}^2}^{p/2-2}
         \norm{Z(s)}^2\trace\CCO{Y(s)QY^*(s)}\,ds
\end{aligned}
\end{multline*}
which yields
\begin{equation}\label{eq:c=zero}
\expect\CCO{\norm{Z(t)}^p}
\leq \frac{p}{2}\CCO{p-1}\expect\int_0^t
        \CCO{\alpha+\norm{Z(s)}^2}^{p/2-1}
          \trace\CCO{Y(s)QY^*(s)}\,ds.
\end{equation}
On the other hand, \eqref{eq:devlop} implies
\begin{multline}\label{eq:cp}
2^{p/2-1}\CCO{c^{p/2}\expect\CCO{V(t)}^{p/2} 
            + \expect\CCO{\alpha+\norm{Z(s)}^2}^{p/2}}\\
\geq \alpha^{p/2}
     +\frac{p}{2}(1+c)\expect\int_0^t \CCO{\alpha+\norm{Z(s)}^2}^{p/2-1}
          \trace\CCO{Y(s)QY^*(s)}\,ds.
\end{multline}
Substituting  \eqref{eq:c=zero} in the right hand side of \eqref{eq:cp}, 
we get that
\begin{equation*}
\frac{\CCO{2c}^{p/2}}{2}\expect\CCO{V(t)}^{p/2}
+2^{p/2-1}\expect\CCO{ \alpha+\norm{Z(t)}^2 }^{p/2}
\geq
{
  \alpha^{p/2}+\frac{1+c}{p-1}\expect\CCO{{\alpha+\norm{Z(t)}^2}}^{p/2} }.
\end{equation*}
Taking the limit when $\alpha$ goes to $0$, 
we then get the result for the case when
$Y$ is uniformly a.e.~bounded. 

In the general case, let us denote, for every integer $N\geq 1$,  
$$Y^N(t)=
\left\{ \begin{array}{ll}
Y(t)                     &\text{ if }\norm{Y(t)}\leq N,\\
N\dfrac{Y(t)}{\norm{Y(t)}}\rule{0em}{1.8em}
                         &\text{ if }\norm{Y(t)}\geq N.
\end{array}\right.$$
By Itô's isometry, $X^N(t):=\int_0^tY^N(s)\,dW(s)$ converges in quadratic
mean to $X(t)$, thus there exists a subsequence 
(still denoted by $(X^N(t))$ for simplicity) 
which converges almost
everywhere to $X$. 
On the other hand, 
if $(e_k)$ is an orthonormal basis of $\h_2$,
we have also, for every $s$, 
\begin{equation*}
\trace\CCO{Y^N(s)Q(Y^N)^*(s)}
=\sum_{k}\norm{Q^{1/2}(Y^N)^*(s)e_k}^2
\nearrow
\trace\CCO{Y(s)QY^*(s)}
\text{ as }N\rightarrow\infty.
\end{equation*}
Using Fatou's and Beppo Levi's lemmas, 
we thus obtain
\begin{align*}
\expect\norm{\int_0^t Y(s)\,dW(s)}^p
= &\expect\CCO{\liminf_N X^N(t)}\\
\leq & \liminf_N\expect\CCO{X^N(t)}\\
\leq & \BDG \liminf_N
              \expect\CCO{ \int_0^t\trace\CCO{Y^N(s)Q(Y^N)^*(s)}\,ds
              }^{p/2}\\
=& \BDG \expect\CCO{ \int_0^t\trace\CCO{Y(s)QY^*(s)}\,ds }^{p/2}.
\end{align*}
\finpr

\begin{prop}\label{cor:Lpboundedness}
With the notations of Theorem \ref{theo:main}  and Lemma \ref{lem:novikov}, 
assume 
that the process $X\in\CUB\bigl(\R,\ellp{2}(\prob, \h_2)\bigl)$ 
satisfies \eqref{eq:mildsol}. Then, for every $p\geq 2$, if
\begin{equation*}
\theta_p':=\frac{2^{3p/2-1}K^{p}}{\delta^{p/2}}
          \CCO{\frac{2^{p/2-1}}{\delta^{p/2}}+\BDG(\trace{Q})^{p/2}}<1,
\end{equation*} 
the family $(X_t)_{t\in\R}$ is bounded in $\ellp{p}$ by a constant
which depends only on $p$, $K$, $\delta$ and $\trace{Q}$.

In particular,  $\theta_2'=\theta'$, 
where $\theta'$ is the constant given in Theorem \ref{theo:main}. 
If $\theta'<1$, we have,  
for every $t\in\R$,
\begin{equation*}
\expect\CCO{\norm{X(t)}^2}
\leq 
\frac{\theta'}{1-\theta'}. 
\end{equation*}
\end{prop}
\proof
We have, using Lemma \ref{lem:novikov}, 
\begin{align*}
\expect\CCO{\norm{X(t)}^p}
\leq &
2^{p-1}\expect{\norm{\int^{t}_{-\infty}S(t-s)F\bigl(s, X(s)\bigl)ds}^p}\\
&+2^{p-1}\expect{\norm{\int^{t}_{-\infty}S(t-s)G\bigl(s,
    X(s)\bigl)dW(s)}^p}\allowbreak\\
\leq &
2^{p-1}\expect\CCO{\int^{t}_{-\infty}e^{-\delta(t-s)}
          \norm{F\bigl(s, X(s)\bigl)}\,ds}^p\\
&+2^{p-1}\BDG\Biggl(\trace{Q}\int^{t}_{-\infty}e^{-2\delta(t-s)}
      \expect { \norm{G\bigl(s, X(s)\bigl)}^2_{L(\h_1, \h_2)}
      }ds\Biggr)^{p/2}\\
\leq &
2^{p-1}\frac{1}{\delta^{p-1}}\expect{\int^{t}_{-\infty}e^{-\delta(t-s)}
          \norm{F\bigl(s, X(s)\bigl)}^p\,ds}\\
&+2^{p-1}\BDG(\trace{Q})^{p/2}\frac{1}{(2\delta)^{p/2-1}}
     \int^{t}_{-\infty}e^{-2\delta(t-s)}
      \expect { \norm{G\bigl(s, X(s)\bigl)}^p_{L(\h_1, \h_2)}
      }ds\\
\intertext{(applying Jensen's inequality under the probabilities
  $\delta e^{-\delta(t-s)}\,ds$ and $2\delta e^{-2\delta(t-s)}\,ds$)}
\leq &
2^{p-1}\frac{1}{\delta^{p-1}}K^p{\int^{t}_{-\infty}e^{-\delta(t-s)}
          \expect(1+\norm{X(s)})^p\,ds}\\
&+2^{p/2}\BDG(\trace{Q})^{p/2}\frac{1}{\delta^{p/2-1}}
     K^p\int^{t}_{-\infty}e^{-2\delta(t-s)}
      \expect(1+\norm{X(s)})^p\,ds\\
\leq & 2^{p-1}\CCO{2^{p-1}\frac{1}{\delta^{p-1}}K^p\int^{t}_{-\infty}e^{-\delta(t-s)}\,ds
+2^{p/2}\BDG(\trace{Q})^{p/2}\frac{1}{\delta^{p/2-1}}
K^p\int^{t}_{-\infty}e^{-2\delta(t-s)}ds}\\
     &+2^{p-1}2^{p-1}\frac{1}{\delta^{p-1}}K^p{\int^{t}_{-\infty}e^{-\delta(t-s)}
          \expect(\norm{X(s)})^p\,ds}\\
&+2^{p-1}2^{p/2}\BDG(\trace{Q})^{p/2}\frac{1}{\delta^{p/2-1}}
     K^p\int^{t}_{-\infty}e^{-2\delta(t-s)}
      \expect(\norm{X(s)})^p\,ds\\
= &
\alpha+
\beta_1\int^{t}_{-\infty}e^{-\delta(t-s)}\expect(\norm{X(s)}^p)\,ds
+\beta_2\int^{t}_{-\infty}e^{-2\delta(t-s)}\expect(\norm{X(s)}^p)\,ds
\end{align*}
with
\begin{gather*}
\alpha
=\frac{2^{3p/2-2}K^{p}}{\delta^{p/2}}
         \CCO{\frac{2^{p/2}}{\delta^{p/2}}+\BDG(\trace{Q})^{p/2}},\\
\beta_1=\frac{2^{2p-2}K^{p}}{\delta^{p-1}},\quad 
\beta_2=\frac{2^{3p/2-1}K^{p}\BDG(\trace{Q})^{p/2}}{\delta^{p/2-1}}.
\end{gather*}
The hypothesis $\theta_p'<1$ is equivalent to $\delta>\beta$, with 
$\beta=\beta_1+\beta_2$. We
conclude by Lemma \ref{lem:gronwall} that 
$$\expect(\norm{X(t)}^p)\leq \alpha\,\frac{\delta}{\delta-\beta}.$$

In the case when $p=2$, 
we have $\bdg{2}=1$, thus 
\begin{equation*}
\alpha=\frac{4K^2}{\delta}
     \CCO{ \frac{1}{\delta} + \frac{\trace{Q}}{2} }
\text{ and }
\beta={4K^2}\CCO{\frac{1}{\delta}+\trace{Q} }
\end{equation*}
and
\begin{equation*}
\expect\CCO{\norm{X(t)}^2}
\leq \frac{ {4K^2}\CCO{ \frac{1}{\delta} + \frac{\trace{Q}}{2} } 
           }{ \delta- {4K^2}\CCO{\frac{1}{\delta}+\trace{Q} }}
\leq \frac{ {4K^2}\CCO{ \frac{1}{\delta} + \trace{Q} } 
           }{ \delta- {4K^2}\CCO{\frac{1}{\delta}+\trace{Q} }}
=\frac{\theta'}{1-\theta'}.
\end{equation*}
\finpr

\begin{rem}\label{rem:novikov-BDG}
We can choose $\BDG$ in Lemma \ref{lem:novikov} such that  
$\lim_{p\rightarrow 2+}\BDG=1$, which implies
$\lim_{p\rightarrow 2+}\theta_p'=\theta'$.
Thus, if $\theta'<1$, and if
$X\in\CUB\bigl(\R,\ellp{2}(\prob, \h_2)\bigl)$ 
satisfies \eqref{eq:mildsol},
the family $(X_t)_{t\in\R}$ is bounded in $\ellp{p}$ for some $p>2$. 
\end{rem}

\proofof{Theorem \ref{theo:main}} Note that 
$$X(t) = \int^{t}_{-\infty}S(t-s)F\bigl(s, X(s)\bigl)ds + 
\int^{t}_{-\infty}S(t-s)G\bigl(s, X(s)\bigl)dW(s)$$ 
satisfies
 $$X(t) = S(t-s)X(s)+ \int^{t}_{s}S(t-s)F\bigl(s, X(s)\bigl)ds 
+ \int^{t}_{s}S(t-s)G\bigl(s, X(s)\bigl)dW(s)$$ 
for all $t\geq s$ for each $s\in \R$ , and hence $X$ is a mild solution
to (\ref{eq:SDE}).

We introduce an operator $L$ by
$$LX(t) = \int^{t}_{-\infty}S(t-s)F\bigl(s, X(s)\bigl)ds + 
\int^{t}_{-\infty}S(t-s)G\bigl(s, X(s)\bigl)dW(s).$$
It can be seen easily that the operator $L$ maps
$\CUB\bigl(\R,\ellp{2}(\prob, \h_2)\bigl)$ into itself.

\medskip
\noindent{\em \textbf{First step.}}
Let us show that $L$ has a unique fixed point. 
We have, for any $t\in\R$, 
\begin{multline*}
 \expect\| (L X)(t) -  (L Y)(t)\|_{\h_2}^2\\
\begin{aligned}
\leq& 2 \expect\Bigl(\int_{-\infty}^{t}e^{-\delta(t-s)}
    \|F(s, X(s))- F(s,Y(s))\|_{\h_2}ds\Bigl)^2\\
  & + 2\expect
\Bigl(\| \int_{-\infty}^{t}S(t-s)[G(s, X(s)) - G(s,Y(s))]dW(s)\|_{\h_2}\Bigl)^2\\
  = &  I_1 + I_2.
\end{aligned}
\end{multline*}
We have
\begin{align*}
I_1  & 
\leq 2\int_{-\infty}^{t}e^{-\delta(t-s)}ds
\int_{-\infty}^{t}e^{-\delta(t-s)}\expect\|F(s, X(s))-F(s,Y(s))\|_{\h_2}^2ds\\
&\leq 2K^2\int_{-\infty}^{t}e^{-\delta(t-s)}ds
       \int_{-\infty}^{t}e^{-\delta(t-s)}\expect\|X(s))- Y(s))\|_{\h_2}^2ds\\
&\leq2K^2\bigl(\int_{-\infty}^{t}e^{-\delta(t-s)}ds\bigl)^2 
          \sup_{s\in\R}\expect\|X(s))- Y(s))\|_{\h_2}^2\\
&\leq\dfrac{2K^2}{\delta^2}\sup_{s\in\R} \expect\|X(s))- Y(s))\|_{\h_2}^2.
\end{align*}
For $I_2$, using the isometry identity we get
\begin{align*}
I_2 &\leq 2\trace Q 
   \int_{-\infty}^{t}e^{-2\delta(t-s)}
         \expect\|G(s, X(s)) - G(s,Y(s))\|^{2}_{L(\h_1, \h_2)}ds\\
&\leq 2\trace QK^2
    \int_{-\infty}^{t}e^{-2\delta(t-s)}
      \expect\|X(s) - Y(s)\|^{2}_{\h_2}ds\\
&\leq 2 K^2\trace Q\bigl(\int_{-\infty}^{t}e^{-2\delta(t-s)}ds\bigl)\,
                        \sup_{s\in\R} \expect\|X(s) - Y(s)\|^{2}_{\h_2}\\
&\leq\dfrac{K^2\trace Q}{\delta}\sup_{s\in\R} \expect\|X(s) - Y(s)\|^{2}_{\h_2}.
\end{align*}
Thus 
$$\expect\|(LX)(t) -  (LY)(t)\|_{\h_2}^2 \leq I_1 + I_2 
           \leq \theta \sup_{s\in\R} \expect\|X(s) - Y(s)\|^{2}_{\h_2}.$$
Consequently, as $\theta < 1$, we deduce that $L$ is a contraction operator, 
hence there exists a unique mild solution to \eqref{eq:SDE} 
in  $\CUB\bigl(\R,\ellp{2}(\prob, \h_1)\bigl)$. 

Furthermore, by \cite[Theorem 7.4]{dapratozabczyk92book}, 
almost all trajectories of 
this solution are 
continuous.

\medskip
\noindent{\em \textbf{Second step}.}
We assume now that $\theta'<1$. 
Let us show that $X$ is almost periodic in distribution. 
We use Bochner's double sequences criterion. 
Let $(\alpha^{'}_n)$ and 
$(\beta^{'}_n)$ be two sequences in $\R$. 
We show that there are
subsequences $(\alpha_n) \subset (\alpha^{'}_n)$ and  
$(\beta_n) \subset (\beta^{'}_n)$ with same indexes such that, for
every $t\in \R$, the limits  
\begin{equation}
 \lim_{n\rightarrow\infty}\lim_{m\rightarrow\infty}\mu(t+\alpha_n+\beta_m)
 \text{ and }\lim_{n\rightarrow\infty}\mu(t+\alpha_n+\beta_n), 
\end{equation}
exist and are equal, where $\mu(t):= \law X(t)$ is the law or
distribution of $X(t)$.

Since $F$ and $G$ are almost periodic, there are subsequences
$(\alpha_n) \subset (\alpha^{'}_n)$ and  
$(\beta_n) \subset (\beta^{'}_n)$ with same indexes such that
\begin{equation}\label{APF}
 \lim_{n\rightarrow\infty}\lim_{m\rightarrow\infty}F(t+\alpha_n+\beta_m, x)
 = \lim_{n\rightarrow\infty}F(t+\alpha_n+\beta_n,x) =: F_0(t,x) 
\end{equation}
and 
\begin{equation}\label{APG}
 \lim_{n\rightarrow\infty}\lim_{m\rightarrow\infty}G(t+\alpha_n+\beta_m,x)
 = \lim_{n\rightarrow\infty}G(t+\alpha_n+\beta_n, x) =: G_0(t,x). 
\end{equation}
These limits exist uniformly with respect to $t \in \R$ and $x$ in
bounded subsets of $\h_2$. 

Set now $(\gamma_n) = (\alpha_n + \beta_n)$. 
For each fixed integer $n$, we consider 
$$X^n(t) = \int^t_{-\infty}S(t-s)F(s+\gamma_n, X^n(s))ds 
     + \int^t_{-\infty}S(t-s)G(s+\gamma_n, X^n(s))dW(s)$$
the mild solution to 
\begin{equation}
 dX^n(t) = AX^n(t)dt + F(t+\gamma_n, X^n(t))dt + G(t+\gamma_n, X^n(t))dW(t)
\end{equation}
and 
$$X^0(t) 
= \int^t_{-\infty}S(t-s)F_0(s, X^0(s))ds 
+ \int^t_{-\infty}S(t-s)G_0(s, X^0(s))dW(s)$$
the mild solution to
\begin{equation}
 dX^0(t) = AX^0(t)dt + F_0(t, X^0(t))dt + G_0(t, X^0(t))dW(t).
\end{equation}
Make the change of variable $\sigma - \gamma_n = s$, the process 
\begin{multline*}
X(t+ \gamma_n) 
= \int^{t+\gamma_n}_{-\infty}S(t+\gamma_n-s)F(s, X(s))ds\\ 
     + \int^{t+\gamma_n}_{-\infty}S(t+\gamma_n-s)G(s, X(s))dW(s)
\end{multline*}
becomes
\begin{multline*}
X(t+ \gamma_n)= \int^t_{-\infty}S(t-s)F(s+\gamma_n, X(s+\gamma_n))ds \\
     + \int^t_{-\infty}S(t-s)G(s+\gamma_n,X(s+\gamma_n))d\tilde{W}_n(s),
\end{multline*} 
where  $\tilde{W}_n(s) = W(s+\gamma_n) - W(\gamma_n)$ 
is a Brownian motion with the same distribution as $W(s)$.
From the independence of the increments of $W$, we deduce that  
the process $X(t+\gamma_n)$ has the same distribution as $X^n(t)$.

Let us show that $X^n(t)$ converges in quadratic 
mean to $X^0(t)$ for each fixed
$t\in \R$. 
We have 
\begin{multline*}
 \expect\lVert X^n(t) - X^0(t)\rVert^2 \\
\begin{aligned}
=& 
\expect\lVert \int^t_{-\infty}S(t-s)[F(s+\gamma_n, X^n(s)) - F_0(s, X^0(s))]ds\\
&+\int^t_{-\infty}S(t-s)[G(s+\gamma_n, X^n(s)) - G_0(s, X^0(s))]dW(s)\rVert^2\\
 \leq&2\expect\lVert \int^t_{-\infty}S(t-s)[F(s+\gamma_n, X^n(s)) 
                      - F_0(s, X^0(s))]ds\rVert^2\\
&+2\expect
\int^t_{-\infty}S(t-s)[G(s+\gamma_n, X^n(s)) - G_0(s, X^0(s))]dW(s)\rVert^2\\
 \leq&4\expect\lVert \int^t_{-\infty}S(t-s)[F(s+\gamma_n, X^n(s)) 
               -F(s+\gamma_n, X^0(s))]ds\rVert^2\\
&+4\expect\lVert\int^t_{-\infty}S(t-s)[F(s+\gamma_n, X^0(s))- F_0(s, X^0(s))]ds\rVert^2\\
&+4\expect\lVert
\int^t_{-\infty}S(t-s)[G(s+\gamma_n, X^n(s)) -G(s+\gamma_n, X^0(s))]dW(s)\rVert^2\\
&+4\expect\lVert
\int^t_{-\infty}S(t-s)[G(s+\gamma_n, X^0(s))-G_0(s, X^0(s))]dW(s)\rVert^2\\
\leq& I_1 + I_2 + I_3+ I_4.
\end{aligned}
\end{multline*} 
Now, using (\ref{cond:semigrp}), (\ref{cond:lipschitz}) 
and the Cauchy-Schwartz inequality, we obtain
\begin{align*}
I_1&= 4\expect\lVert \int^t_{-\infty}S(t-s)
      [F(s+\gamma_n, X^n(s)) -F(s+\gamma_n, X^0(s))]ds\rVert^2\\
&\leq4\expect\bigl(\int^t_{-\infty}\lVert 
  S(t-s)\rVert \lVert F(s+\gamma_n, X^n(s)) -F(s+\gamma_n, X^0(s))
  \rVert ds\bigl)^2\\
&\leq4\expect\bigl(\int^t_{-\infty}e^{-\delta(t-s)} 
  \lVert F(s+\gamma_n, X^n(s)) -F(s+\gamma_n, X^0(s))\rVert ds\bigl)^2\\
&\leq4\bigl(\int^t_{-\infty}e^{-\delta(t-s)}ds\bigl)
    \bigl(\int^t_{-\infty}e^{-\delta(t-s)}\expect\lVert F(s+\gamma_n, X^n(s))
               -F(s+\gamma_n, X^0(s))\rVert^2 ds\bigl)\\
&\leq \frac{4K^2}{\delta}\int^t_{-\infty}e^{-\delta(t-s)} 
                                \expect\lVert X^n(s) - X^0(s)\rVert^2ds.
\end{align*}
Then we obtain
\begin{align*}
 I_2 & = 4\expect\lVert\int^t_{-\infty}S(t-s)[F(s+\gamma_n, X^0(s))
          - F_0(s, X^0(s))]ds\rVert^2\\
&\leq 4\expect\bigl(\int^t_{-\infty}e^{-\delta(t-s)} 
     \lVert F(s+\gamma_n, X^0(s)) -F_0(s, X^0(s))\rVert ds\bigl)^2\\
&\leq 4\expect\bigl(\int^t_{-\infty}e^{-\delta(t-s)}ds\bigl)
\bigl(\int^t_{-\infty}e^{-\delta(t-s)} 
         \lVert F(s+\gamma_n, X^0(s)) -F_0(s, X^0(s))\rVert^2
  ds\bigl)\\
&=\frac{4}{\delta}
   \expect\bigl(\int^t_{-\infty}e^{-\delta(t-s)} 
         \lVert F(s+\gamma_n, X^0(s)) -F_0(s, X^0(s))\rVert^2
  ds\bigr)\\
&=
  2
  \expect\bigl(\int^t_{-\infty}\frac{2e^{-\frac{\delta}{2}(t-s)}}{\delta} 
     \Bigl(e^{-\frac{\delta}{2}(t-s)} 
         \lVert F(s+\gamma_n, X^0(s)) -F_0(s, X^0(s))\rVert^2\Bigr)
                                                                 ds\bigr).
\end{align*}
Since $X^0\in\CUB\Bigl(\R,\ellp{2}(\prob, \h_2)\Bigr)$ and 
$\sup_{t\in\R} \expect\lVert X^0(t)\rVert^2 < \infty$, 
the family
$$\Bigl(e^{-\frac{\delta}{2}(t-s)}\norm{X^0(s)}^2\Bigr)_{-\infty<s\leq t}$$
is uniformly integrable. 
Indeed, for any sequence $(s'_n)$ in $(-\infty,t]$, there exists a
subsequence $(s_n)$ which converges to some $s\in[-\infty,t]$. If
$s>-\infty$, the sequence $\bigl(e^{-\frac{\delta}{2}(t-s_n)}X^0(s_n)\bigr)$ 
converges in $\ellp{2}(\prob, \h_2)$ to $e^{-\frac{\delta}{2}(t-s)}X^0(s)$,
and if $s=-\infty$, it converges to $0$. 
Thus any sequence $\bigl(e^{-\frac{\delta}{2}(t-s_n)}X^0(s_n)\bigr)$
contains a subsequence which is convergent in $\ellp{2}(\prob, \h_2)$,
which proves the uniform integrability. 
Alternatively, one can use Remark \ref{rem:novikov-BDG}, 
since $\theta'<1$, which yields 
uniform integrability of 
$\bigl(\norm{X^0(t)}^2\bigr)$.
By the growth condition \eqref{cond:croissance}, this shows that the
family 
$$(U_{s,n}):=\Bigl(e^{-\frac{\delta}{2}(t-s)} 
         \lVert F(s+\gamma_n, X^0(s)) -F_0(s, X^0(s))\rVert^2
                                                  \Bigr)_{-\infty<s\leq
                                                    t,\, n\geq 1}$$
is uniformly integrable. By La Vall\'ee Poussin's criterion, there
exists a non-negative increasing convex function $\Phi :
\R\rightarrow\R$ 
such that
$\lim_{t\rightarrow\infty}\frac{\Phi(t)}{t}=+\infty$ and
$\sup_{s,n}\expect(\Phi(U_{s,n}))<+\infty$. 
We thus have 
$$\sup_{n}\expect \int^t_{-\infty}\frac{2e^{-\frac{\delta}{2}(t-s)}}{\delta} 
     \Phi\bigl(U_{s,n}\bigr)\,ds <+\infty,$$
which prove that the family
$(U_{.,n})_{n\geq 1}$
is uniformly integrable. with respect to the probability measure
$\prob\otimes \frac{2}{\delta}e^{-\frac{\delta}{2}(t-s)}ds$ on
$\esprob\times(-\infty,t]$.  
We deduce by \eqref{APF} 
that $I_2$ converges to $0$ as $n\rightarrow\infty$.

Applying Itô's isometry, we get
\begin{align*}
 I_3& = 4\expect\lVert\int^t_{-\infty}S(t-s)
    [G(s+\gamma_n, X^n(s)) -G(s+\gamma_n, X^0(s))]dW(s)\rVert^2\\
& \leq 4\trace Q\expect\int^t_{-\infty}\lVert S(t-s)\rVert^2 
      \lVert G(s+\gamma_n, X^n(s)) -G(s+\gamma_n, X^0(s))\rVert^2 ds\\
&\leq \frac{4}{\delta}\trace Q\int^t_{-\infty}e^{-2\delta(t-s)} \expect\lVert 
        G(s+\gamma_n, X^n(s)) -G(s+\gamma_n, X^0(s))\rVert^2 ds\\
&\leq 4K^2\trace Q\int^t_{-\infty}e^{-2\delta(t-s)} 
      \expect\lVert X^n(s) - X^0(s)\rVert^2ds.
\end{align*}
and
\begin{align*}
 I_4& = 4\expect\lVert\int^t_{-\infty}S(t-s)
     [G(s+\gamma_n, X^0(s))-G_0(s, X^0(s))]dW(s)\rVert^2\\
& \leq 4\trace Q\expect\bigl(\int^t_{-\infty}\lVert S(t-s)\rVert^2 
      \lVert G(s+\gamma_n, X^0(s))-G_0(s, X^0(s))\rVert^2 ds\bigl)\\
& \leq 4\trace Q\expect\bigl(\int^t_{-\infty}e^{-2\delta(t-s)}
      \lVert G(s+\gamma_n, X^0(s))-G_0(s, X^0(s))\rVert^2 ds\bigl).
\end{align*}
For the same reason as for $I_2$, 
the right hand term goes to $0$ as
$n\rightarrow\infty$. 

We thus have 
\begin{multline*}
\expect\lVert X^n(t) - X^0(t)\rVert^2
\leq \alpha_n 
  +\frac{4K^2}{\delta}\int^t_{-\infty}e^{-\delta(t-s)} \expect\lVert
  X^n(s)-X^0(s)\rVert^2\,ds\\
  +4K^2\trace{Q}\int^t_{-\infty}e^{-2\delta(t-s)} \expect\lVert
  X^n(s)-X^0(s)\rVert^2\,ds
\end{multline*}
for a sequence $(\alpha_n)$ such that
$\lim_{n\rightarrow\infty}\alpha_n=0.$ Furthermore,  
$\beta:=\frac{4K^2}{\delta}+4K^2\trace{Q}<\delta$. We conclude by 
 Lemma \ref{lem:gronwall} that 
$$\lim_{n\rightarrow\infty}  \expect\lVert X^n(t) - X^0(t)\rVert^2 = 0,$$
hence $X^n(t)$ converges in distribution to $X^0(t)$. 
But, since the distribution of $X^n(t)$ is the same as that of $X(t+\gamma_n)$,
 we deduce that $X(t+\gamma_n)$ converges in distribution to  $X^0(t)$, i.e. 
$$ \lim_{n\rightarrow\infty}\mu(t+\alpha_n+\beta_n) =\law{X^0(t)}
   =: \mu^0_t .$$
By analogy and using (\ref{APF}), (\ref{APG}) we can easily deduce that
 $$\lim_{n\rightarrow\infty}\lim_{m\rightarrow\infty}\mu(t+\alpha_n+\beta_m) = \mu^0_t.$$

We have thus proved that $X$ 
has almost periodic one-dimensional distributions. 
To prove that $X$ is almost periodic in distribution, we apply 
Proposition \ref{prop:daprato-tudor3.1}: for fixed $\tau\in\R$, let   
$\xi_n=X(\tau+\alpha_n)$, $F_n(t,x)=F(t+\alpha_n,x)$,
$G_n(t,x)=G(t+\alpha_n,x)$. 
By the foregoing, $(\xi_n)$ converges in distribution to some variable
$Y(\tau)$. 
We deduce that $(\xi_n)$ is tight, 
and thus $(\xi_n,W)$ is tight also. 
We can thus choose $Y(\tau)$ such that $(\xi_n,W)$ converges in
distribution to $(Y(\tau),W)$. 
Then, by Proposition \ref{prop:daprato-tudor3.1}, for every $T\geq \tau$, 
$X(.+\alpha_n)$ converges in distribution on $C([\tau,T]; \h_2)$ to
the (unique in distribution) solution to 
\begin{equation*}
Y(t)=S(t-\tau)Y(\tau)+\int^{t}_{\tau}S(t-s)F\bigl(s, Y(s)\bigl)\,ds 
+ \int^{t}_{\tau}S(t-s)G\bigl(s, Y(s)\bigl)\,dW(s).
\end{equation*}
Note that $Y$ does not depend on the chosen interval $[\tau,T]$, thus
the convergence takes place on $C(\R; \h_2)$. 
Similarly, $Y_n:=Y(.+\beta_n)$ converges in distribution on
$C(\R; \h_2)$ to a continuous process $Z$ such that, for $t\geq\tau$, 
$$Z(t)=S(t-\tau)Z(\tau)+\int^{t}_{\tau}S(t-s)F\bigl(s, Z(s)\bigl)ds + 
\int^{t}_{\tau}S(t-s)G\bigl(s, Z(s)\bigl)dW(s).$$
But, by \eqref{APF} and \eqref{APG}, $X(.+\gamma_n)$ converges in
distribution to the same process $Z$. 
Thus $X$ is almost periodic in distribution. 
\finpr

\section{Weak averaging}\label{sect:averaging}
In this section, we strengthen slightly the assumptions on the
semigroup $S$ but we replace the condition of almost periodicity on $F$
and $G$ by a weaker condition that keeps only the features of almost
periodicity that are useful for averaging. 
More precisely, 
we assume Conditions 
\eqref{cond:wienerpro},
\eqref{cond:croissance}, and \eqref{cond:lipschitz}
of Section \ref{sect:Form.prob}, 
but we replace Condition \eqref{cond:semigrp} by the stronger 
Condition (\ref{cond:semigrp-}') below and 
Condition \eqref{cond:pp}
 by the weaker Condition (\ref{cond:pp}') below :
\begin{enumerate}[{\rm (i')}]
\setcounter{enumi}{\value{carbon2}}
  \item \label{cond:semigrp-} Condition \eqref{cond:semigrp} is satisfied
    and the semigroup $S$ is {immediately norm continuous} 
(see \cite[Definition II.4.24]{engel-nagel00book}), 
i.e.~the mapping $t\mapsto S(t)$ is continuous in
    operator norm on $]0,\infty]$.  
\end{enumerate}
\begin{enumerate}[{\rm (i')}]
\setcounter{enumi}{\value{carbone}}
  \item \label{cond:pp-}
The mappings 
$F:\R \times \h_2\rightarrow \h_2$ and 
$G:\R\times\h_2\rightarrow L(\h_1, \h_2)$ satisfy :
\begin{enumerate}[(a)]  
\item\label{imagecompacte} For every compact subset $K$ of $\h_2$, the sets
$$\accol{F(t,x)\tq t\in\R,\ x\in K}
\text{ and }
\accol{G(t,x)\tq t\in\R,\ x\in K}$$
are compact.  
\item There exist continuous functions $F_0: \h_2 \rightarrow \h_2$  and 
$G_0:\h_2\rightarrow L(\h_0, \h_2)$ 
satisfying \eqref{eq:moyenne_simple} and \eqref{eq:gaussassumption} 
uniformly with respect to $x$ in compact subsets of $\h_2$.
\end{enumerate}\end{enumerate}
Contrarily to \cite{vrkoc95weakaveraging}, no condition of analyticity
of $S$ is required. Condition (\ref{cond:semigrp-}') is satisfied by a
broad class of semigroups, see 
\cite{engel-nagel00book} for details.
Condition (\ref{cond:pp}') 
is weaker than \eqref{cond:pp}
thanks to Propositions \ref{prop: mean} and \ref{rem:mean}.

Let us define the Hilbert space
$\h_0=\range{Q^{1/2}}$ as in Proposition \ref{rem:mean}, where $Q$ is the
covariance operator of the Wiener process $W$.

\begin{lem}\label{lem:moyenne_S}
Under Hypothesis 
\eqref{cond:wienerpro}, 
{\rm(\ref{cond:semigrp-}')}, 
\eqref{cond:croissance},
\eqref{cond:lipschitz}, and 
{\rm(\ref{cond:pp-}')}, for any continuous function
$x :\,\R\rightarrow\h_2$, we have 
\begin{align}
&\lim_{\epsilon\rightarrow 0^+}\int_{\Tau}^{\Tau+t} S(\Tau+t-s)
F\Bigl(\frac{s}{\epsilon},x(s)\Bigr)\,ds
=\int_{\Tau}^{\Tau+t} S(\Tau+t-s) F_0(x(s))\,ds\label{eq:SmoyenneF}\\
& \lim_{\epsilon\rightarrow 0^+}\Bigl|
\int^{\Tau+t}_{\Tau}\Bigl\lgroup
S(\Tau+t-s) G\Bigl(\frac{s}{\epsilon},x(s)\Bigr)
   QG^*\Bigl(\frac{s}{\epsilon},x(s)\Bigr) S^*(\Tau+t-s)\label{eq:SmoyenneG}\\
&\phantom{ \lim_{\epsilon\rightarrow 0^+}\Bigl|\int^{\Tau+t}_{\Tau}\Bigl( }           
     - S(\Tau+t-s)G_0(x(s))QG^{*}_{0}(x(s))S^*(\Tau+t-s)\Bigr\rgroup\,ds
                          \Bigr|_{\NN} = 0\notag
\end{align}
for all $\Tau\in\R$ and $t>0$.
\end{lem}
\proof
Let $\Tau\in\R$ and $t>0$. 
Let $\gamma>0$, and let us choose
$\alpha>0$ such that 
\begin{align*}
&\norm{\int_{\Tau+t-\alpha}^{\Tau+t}
             S(\Tau+t-s)F\Bigl(\frac{s}{\epsilon},x(s)\Bigr)\,ds
}<\gamma
\intertext{and}
&\norm{\int_{\Tau+t-\alpha}^{\Tau+t}
             S(\Tau+t-s)F_0(x(s))\,ds
}<\gamma.
\end{align*}
This is possible since the functions inside the integrals are bounded. 
By Condition (\ref{cond:semigrp-}'), the semigroup $S$ is
uniformly continuous on $[\alpha,t]$. 
We can thus divide the interval $[\Tau,\Tau+t-\alpha]$ by a partition 
$(\Tau+ih)_{i=0,\dots,N}$, in such a way that 
$\norm{S(\Tau+t-s)-S(t+ih)}<\gamma$ for $s\in[\Tau+ih,\Tau+(i+1)h]$, 
$i=0,\dots,N-1$. We have
\begin{align*}
\int_{\Tau}^{\Tau+t-\alpha}
             S(\Tau+t-s)F\Bigl(\frac{s}{\epsilon},x(s)\Bigr)\,ds
&=\sum_{i=0}^{N-1}
\int_{\Tau+ih}^{\Tau+(i+1)h}
             S(\Tau+t-s)F\Bigl(\frac{s}{\epsilon},x(s)\Bigr)\,ds
\intertext{and}
\int_{\Tau}^{\Tau+t-\alpha}
             S(\Tau+t-s)F_0(x(s))\,ds
&=\sum_{i=0}^{N-1}
\int_{\Tau+ih}^{\Tau+(i+1)h}
             S(\Tau+t-s)F_0(x(s))\,ds.
\end{align*}
But, for some constant $C$ and $i=0,\dots,N-1$, 
\begin{align*}
&\norm{\int_{\Tau+ih}^{\Tau+(i+1)h}
             S(\Tau+t-s)F\Bigl(\frac{s}{\epsilon},x(s)\Bigr)\,ds
-\int_{\Tau+ih}^{\Tau+(i+1)h}
             S(t-ih)F\Bigl(\frac{s}{\epsilon},x(s)\Bigr)\,ds}
\leq C\gamma
\intertext{and}
&\norm{\int_{\Tau+ih}^{\Tau+(i+1)h}
             S(\Tau+t-s)F_0(x(s))\,ds
-\int_{\Tau+ih}^{\Tau+(i+1)h}
             S(t-ih)F_0(x(s))\,ds}
\leq C\gamma.
\end{align*}
Now, by the Krasnoselski-Krein 
lemma \cite{krasnoselski-krein56}, 
we have
\begin{equation*}
\lim_{\epsilon\rightarrow 0}
S(t-ih)\int_{\Tau+ih}^{\Tau+(i+1)h}F\Bigl(\frac{s}{\epsilon},x(s)\Bigr)\,ds
=S(t-ih)\int_{\Tau+ih}^{\Tau+(i+1)h}F_0(x(s))\,ds,
\end{equation*}
which proves \eqref{eq:SmoyenneF}.

To prove \eqref{eq:SmoyenneG}, we need to show that 
\begin{align*}
\lim_{\epsilon\rightarrow 0}
\Biggl(
&\int_{\Tau}^{\Tau+t}
\Bigl\langle S(\Tau+t-s)G\Bigl(\frac{s}{\epsilon},x(s)\Bigr)
   QG^*\Bigl(\frac{s}{\epsilon},x(s)\Bigr)
   S^*(\Tau+t-s)x,y\Bigr\rangle\,ds\\
&-\int_{\Tau}^{\Tau+t}
\Bigl\langle S(\Tau+t-s)G_0(x(s))
   QG_0^*(x(s))
   S^*(\Tau+t-s)x,y\Bigr\rangle\,ds
\biggr)=0.
\end{align*}
uniformly with respect to $x,y$ in the unit ball $B(0,1)$ of $\h_2$. 
As in the proof of \eqref{eq:SmoyenneF}, we have
\begin{align*}
&\int_{\Tau}^{\Tau+t-\alpha}
\Bigl\langle G\Bigl(\frac{s}{\epsilon},x(s)\Bigr)
   QG^*\Bigl(\frac{s}{\epsilon},x(s)\Bigr)
   S^*(\Tau+t-s)x,S^*(\Tau+t-s)y\Bigr\rangle\,ds\\
&=\sum_{i=0}^{N-1}\int_{\Tau+ih}^{\Tau+(i+1)h}
\Bigl\langle G\Bigl(\frac{s}{\epsilon},x(s)\Bigr)
   QG^*\Bigl(\frac{s}{\epsilon},x(s)\Bigr)
   S^*(\Tau+t-s)x,S^*(\Tau+t-s)y\Bigr\rangle\,ds
\intertext{and}
&\int_{\Tau}^{\Tau+t-\alpha}
\Bigl\langle G_0(x(s))
   QG_0^*(x(s))
   S^*(\Tau+t-s)x,S^*(\Tau+t-s)y\Bigr\rangle\,ds\\
&=\sum_{i=0}^{N-1}\int_{\Tau+ih}^{\Tau+(i+1)h}
\Bigl\langle G_0(x(s))
   QG_0^*(x(s))
   S^*(\Tau+t-s)x,S^*(\Tau+t-s)y\Bigr\rangle\,ds.
\end{align*}
Replacing in the right hand sides $S^*(\Tau+t-s)$ by $S^*(t-ih)$,
we reduce the proof to the equation
\begin{align*}
 \lim_{\epsilon\rightarrow 0}
&\int_{\Tau+ih}^{\Tau+(i+1)h}
\Bigl\langle G\Bigl(\frac{s}{\epsilon},x(s)\Bigr)
   QG^*\Bigl(\frac{s}{\epsilon},x(s)\Bigr)
   S^*(t-ih)x,S^*(t-ih)y\Bigr\rangle\,ds\\
&-\int_{\Tau+ih}^{\Tau+(i+1)h}
\Bigl\langle G_0(x(s))
   QG_0^*(x(s))
   S^*(t-ih)x,S^*(t-ih)y\Bigr\rangle\,ds=0
\end{align*}
uniformly with respect to $x,y$ in $B(0,1)$. But again this follows
from the Krasnoselski-Krein lemma.
\finpr

Recall that, 
if $X$ and $Y$ are two random vectors of a Banach space $\Bspace$, 
the $\ellp{2}$-Wasserstein distance $\WASS(X,Y)$ 
between the distributions of $X$ and $Y$ is 
$$\WASS(X,Y)
=\left(
\inf
       \expect\CCO{\norm{\widehat{X}-\widehat{Y}}^2_{\Bspace}}
\right)^{1/2}$$
where the infimum is taken over all joint distributions of random
vectors $\widehat{X}$ and $\widehat{Y}$ satisfying 
$\law{\widehat{X}}=\law{X}$ and $\law{\widehat{Y}}=\law{Y}$.

By e.g.~\cite[Theorem 6.9]{villani09oldnew}), 
if $(X_n)$ is a sequence of random vectors of $\Bspace$ and if $X$ is a
random vector of $\Bspace$, 
the sequence $(\law{X_n})$ converges to $\law{X}$ for $\WASS$
if and only if 
$(X_n)$ 
converges to $X$ in distribution and 
$\CCO{\norm{X_n}^2_{\Bspace}}$ is uniformly
integrable.

If $X$ and $Y$ are continuous $\h_2$-valued stochastic processes, 
for any interval $[a,b]$, we denote by $\WASS_{[a,b]}$ 
the $\ellp{2}$-Wasserstein distance between the ditributions of $X$
and $Y$, seen as 
$\C([a,b], \h_2)$-valued random
variables.

We are now ready to state our main averaging result.

\begin{theo}\label{theo:averaging}
Let the assumptions 
\eqref{cond:wienerpro}, 
{\rm(\ref{cond:semigrp-}')}, 
\eqref{cond:croissance},
\eqref{cond:lipschitz}, 
 and  
{\rm(\ref{cond:pp}')} 
be fulfilled 
and the constant 
 $\theta' =\dfrac{4K^2}{\delta}\CCO{\dfrac{1}{\delta} + \trace Q}
<1$. 
For each fixed $\varepsilon \in ]0, 1[$,
 let $X^\varepsilon$ be the mild solution to the equation
\begin{equation}\label{eq:ASDE}
 dX^\varepsilon(t) 
= AX^\varepsilon(t)dt + F\Bigl(\frac{t}{\varepsilon}, X^\varepsilon(t)\Bigr)\,dt + 
      G\Bigl(\frac{t}{\varepsilon}, X^\varepsilon(t)\Bigr)\,dW(t),
\end{equation} 
and let
$X^{0}$ be the mild solution to
\begin{equation}\label{eq:station}
 d X^{0}(t) = A(X^{0}(t))dt + F_0(X^{0}(t))dt + G_0(X^{0}(t))dW(t),
 \end{equation}
which is a stationary process.
Then $\WASS_{[a,b]}\CCO{X^{\varepsilon},X^{0}}\rightarrow 0$ 
as 
$\varepsilon\rightarrow 0+$,
for any compact interval $[a,b]$. 
\end{theo}
Before we give the proof of this theorem, let us recall some well-known results.

\begin{prop}\label{prop: Gaussian}(\cite{Chevetcompact83})
 Let $(X_n)_{n\geq0}$ be a sequence of centered Gaussian random
 variable on a separable Hilbert space $\h$  
with sequence of covariance operators $(Q_n)_{n\geq0}$. Then
$(X_n)_{n\geq0}$ converges in distribution to $X_0$ in $\h$  
if and only if $$\mid Q_n - Q_0\mid_\NN \rightarrow0, n\rightarrow\infty$$
\end{prop}

Let $\polish, \V, \h$ be real separable Hilbert spaces, let $W$ be a
$\polish$-valued $(\tribu_t)$-adapted Wiener 
 process with nuclear covariance operator $Q$.
\begin{prop}\label{prop: Randomvar}
(\cite[Proposition 2.2]{vrkoc95weakaveraging})
Let $\alpha: \h\rightarrow \V$ be a Lipschitz mapping and 
$\sigma: \R \times \h \rightarrow L(\polish, \V)$ 
be a measurable mapping such that 
$\|\sigma(r,x)\|_{L(\polish, \V)} \leq M (1 + \|x \|_{\h})$ and  
$\|\sigma(r,x) - \sigma(r,y)\|_{L(\polish, \V)} \leq M \|x - y\|_{\h}$ 
for a constant $M$ and every 
$ r \in  [s,t], x, y \in \h $. 
Let $g\in \BL(\V)$, we define 
$$\psi(y) = \expect g\Bigl(\alpha(y) 
           + \int^{t}_{s}\sigma(r,y)dW(r)\Bigl), y\in \h .$$ 
Let $u: \esprob\rightarrow\h$ be a $(\tribu_s)$-measurable 
random variable with $\expect\|u\|^{2}_{\h}< \infty$. 
Then 
$$\expect\Bigl[g\Bigl(\alpha(u) 
      + \int^{t}_{s}\sigma(r,u)dW(r)\Bigl)|\tribu_s\Bigl] 
= \psi(u)\, \,  \prob\text{-a.s.}$$ 
\end{prop}

\proofof{Theorem \ref{theo:averaging}}
We denote 
$F_\varepsilon(s,x):= F\Bigl(\dfrac{s}{\varepsilon}, x\Bigr)$,
 $G_\varepsilon(s,x):= G\Bigl(\dfrac{s}{\varepsilon}, x\Bigr)$, 
and, for every $X\in \CUB\Bigl(\R,\ellp{2}(\prob, \h_2)\Bigr) $, 
\begin{gather*}
L_\varepsilon(X)(t)
:=\int_{-\infty}^{t}S(t-s)F_\varepsilon(s,X(s))ds
          + \int_{-\infty}^{t}S(t-s)G_\varepsilon(s,X(s))dW(s),\\
L_0(X)(t):=\int_{-\infty}^{t}S(t-s)F_0(X(s))ds
          +\int_{-\infty}^{t}S(t-s)G_0(X(s))dW(s).
\end{gather*}

\medskip
\noindent{\em \textbf{First step.}}
Let $X \in \CUB\bigl(\R,\ellp{2}(\prob, \h_2)\bigl)$. 
Let us show that $L_\varepsilon(X)\rightarrow L_0(X)$ in distribution, as
$\varepsilon\rightarrow0$,  
in the space $\C(\R, \h_2)$ endowed with the topology of uniform
convergence on compact intervals of $\R$. 
This amounts to prove that, for any compact interval $[a,b]$,  
 $L_\varepsilon(X)\rightarrow L_0(X)$ in distribution 
in the space $\C([a,b], \h_2)$  (see \cite[Theorem 5]{Whitt70}).
Actually, we will prove a slightly stronger result, using the
$\ellp{2}$-Wasserstein distance.

By Conditions \eqref{cond:semigrp} and \eqref{cond:croissance} 
we deduce that, for every $\eta > 0$, there
exists $\Tau$  such that, for every $\eta\geq 0$ and for each $\sig<\Tau$,  
$\expect\|L_\varepsilon(X)(\sig)\|^{2}_{\h_2}< \eta$.
Thus, for the proof that
$L_\varepsilon(X)$ converges in distribution to $L_0(X)$ on $\C(\R, H_2)$, 
it suffices to show the convergence in distribution on 
$\C([\Tau, T], \h_2)$, for every $T\geq\Tau$, of 
$$Y^\epsilon(t):=\int_{\Tau}^{t}S(t-s)F_\varepsilon(s,X(s))\,ds
          + \int_{\Tau}^{t}S(t-s)G_\varepsilon(s,X(s))\,dW(s)$$ 
to
$$
Y^0(t):=\int_{\Tau}^{t}S(t-s)F_0(X(s))\,ds
          +\int_{\Tau}^{t}S(t-s)G_0(X(s))\,dW(s).
$$

As  $X \in \CUB\bigl(\R,\ellp{2}(\prob, \h_2)\bigl)$, 
it satisfies the following condition: 
For every $\eta>0$, there exist a partition 
$$\{\Tau = t_o<t_1< \dots<t_k = T\}
\text{ of }  [\Tau, T]$$ 
and an adapted process $$\transl{X}(t) = \Sigma_{i = 0}^{k-1}
{X}(t_{i})1_{[t_{i1}, t_{i+1}[}(t)$$
such that 
$$\sup_{t\in[\Tau,T]} \expect\|X(t) - \transl{X}(t)\|^{2}_{\h_2}< \eta.$$ 
Using the fact that $L_\varepsilon$ is Lipschitz, we can furthermore
choose the partition  $(t_o,\dots,t_k)$ such that
\begin{equation}\label{eq:partition_uniforme}
\sup_{\epsilon>0}\sup_{t\in[\Tau,T]} 
\expect\|L_\varepsilon X(t) - L_\varepsilon
\transl{X}(t)\|^{2}_{\h_2}<\eta.
 \end{equation}

For $\epsilon>0$, we denote 
$F_\epsilon(s,x):=F(\frac{s}{\varepsilon},x)$, 
$G_\epsilon(s,x):= G(\frac{s}{\varepsilon},x)$,
and we set 
\begin{align*}
 \transl{X}^\epsilon(t) 
=&\int_\Tau^t S(t-s)F_\epsilon(s,\transl{X}({s}))\,ds
+\int_\Tau^t S(t-s)G_\epsilon(s,\transl{X}({s}))\,dW(s)
\\
=&\sum_{i = 0}^{k-1} \Bigl(\int_{t_{i}\wedge t}^{t_{i+1}\wedge t}S(t-s)F_\epsilon(s,
{X}({t_{i}}))\,ds\\
& +\int_{t_{i}\wedge t}^{t_{i+1}\wedge t}S(t-s)G_\epsilon(s,{X}({t_{i}}))\,dW(s)\Bigl).
\end{align*}
By \cite[Theorem 6.10]{dapratozabczyk92book}, each $\transl{X}^\epsilon$ has
a continuous modification. 

Let us prove that, for each $l=1,\dots,k$,  
$\transl{X}_{t_1}^{\epsilon}$ converges in distribution to 
$\transl{X}_{t_1}^{0}$
as $\epsilon\rightarrow 0$.
We define a mapping  
$$
\gamma_\epsilon : \h_{2}^{l}\rightarrow  \ellp{1}(\esprob, \h_2)
$$
by 
 \begin{multline*}
  \gamma_\epsilon(y_0, y_1, \dots, y_{l-1})
=\sum_{i = 1}^{l}
  \Bigl(\int_{t_{i-1}}^{t_{i}}S(t_l-s)F_\epsilon(s,y_{i-1})\,ds\\ 
+ \int_{t_{i-1}}^{t_{i}}S(t_l-s)G_\epsilon(s,y_{i-1})\,dW(s)\Bigl).
 \end{multline*}
Using Proposition \ref{prop: Randomvar}, we get that 
$$\law{\transl{X}_{t_l}^{\epsilon}}
=\law{\gamma_\epsilon(\transl{X}_{t_0}, \transl{X}_{t_1}, 
                        \dots,\transl{X}_{t_{l-1}})}.$$
Let 
$$\mu_{t_0,t_1, \dots,t_{l-1} }
=\law{\transl{X}_{t_0}, \transl{X}_{t_1}, \dots,
  \transl{X}_{t_{l-1}}}.$$
Let $g\in\BL(\h_2)$, 
and $h_\epsilon(y) = \expect[g(\gamma_\epsilon(y))];\,\, y\in  \h_{2}^{l}$.
We have 
$$\expect[g(\transl{X}_{t_l}^{\epsilon})] = \expect[h_\epsilon(\transl{X}_{t_0},
\transl{X}_{t_1}, \dots, \transl{X}_{t_{l-1}})] = \int_{\h_{2}^{l}}
h_\epsilon(y)d\mu_{t_0,t_1, \dots,t_{l-1} }(y)$$
and 
 \begin{align*}
  |\expect[g(\transl{X}_{t_l}^{\epsilon})] - \expect[g(\transl{X}_{t_l}^{0})]|
& = |\int_{\h_{2}^{l}}  h_\epsilon(y)- h_0(y)d\mu_{t_0,t_1, \dots,t_{l-1} }(y)|\\
&\leq\int_{\h_{2}^{l}}  |h_\epsilon(y)- h_0(y)|d\mu_{t_0,t_1, \dots,t_{l-1} }(y).
 \end{align*}
Let us show that 
$h_\epsilon(y)\rightarrow h_0(y)$ as $\epsilon\rightarrow 0$ 
for every $y\in\h_{2}^{l}$. We have   
 \begin{align*}
 \gamma_\epsilon(y)- \gamma_0(y) 
= &\sum_{i = 1}^{l} \int_{t_{i-1}}^{t_{i}}S(t_l-s)
                \bigl(F_\epsilon(s,y_{i-1})- F_0(y_{i-1})\bigl)ds\\
& + \sum_{i = 1}^{l}
 \int_{t_{i-1}}^{t_{i}}S(t_l-s)\bigl(G_\epsilon(s,y_{i-1})- G_0(y_{i-1})\bigl)\,dW(s)\\
=& I_{\epsilon} + J_{\epsilon}.
 \end{align*}
Lemma \ref{lem:moyenne_S} 
implies  that 
$I_{\epsilon}\rightarrow0$ as $\epsilon\rightarrow0$, and since
 $$\sum_{i = 1}^{l}\int_{t_{i-1}}^{t_{i}}S(t_l-s)G_\epsilon(s,y_{i-1})\,dW(s)$$ 
is a centered Gaussian random variable in
 $\h_{2}$, we deduce by 
Lemma \ref{lem:moyenne_S} 
and Proposition \ref{prop: Gaussian}  
that $J_{\epsilon}\rightarrow0$ 
in distribution as $\epsilon\rightarrow0$ 
hence   $\gamma_\epsilon(y)\rightarrow \gamma_0(y)$ in distribution 
as $\epsilon\rightarrow0$. 
Consequently
  \begin{equation}\label{eq:converj}
    h_\epsilon(y)\rightarrow h_0(y)\text{ for any }  y\in\h_{2}^{l}. 
  \end{equation}

For every $\eta^{'}>0$, there exists a compact set $\FF \subset
\h_{2}^{l}$ such that 
 $$\mu_{t_0,t_1, \dots,t_{l-1} }(\h_{2}^{l}\setminus \FF)<\eta^{'}.$$ 
We have 
\begin{equation}\label{eq:Lipschitz}
 h_\epsilon \in \BL(\h_{2}^{l}) \text{ and }  
\sup_\epsilon\|h_\epsilon\|_{\bl} <\infty
\end{equation}
because, for all $y, z \in \h_{2}^{l}$, and for some constant $K_1$,
\begin{align*}
 |h_\epsilon(y)- h_\epsilon(z)| 
 \leq \|g\|_{\bl}\expect\| \gamma_\epsilon(y)- \gamma_\epsilon(z)\|_{\h_{2}} 
 \leq K_1\|g\|_{\bl}\|y-z\|_{\h_{2}^{l}}.
\end{align*}
From (\ref{eq:converj}), (\ref{eq:Lipschitz}) and the compactness of
$\FF$, 
we deduce that $h_\epsilon$ converges to $h_0$ 
uniformly 
on $\FF$, hence 
 $$\lim_{\epsilon\rightarrow0}
\int_{\FF}  |h_\epsilon(y)- h_0(y)|d\mu_{t_0,t_1, \dots,t_{l-1} }(y) = 0$$  
and, since $g$ is a bounded function,
\begin{multline*}
\int_{\h_{2}^{l}\setminus\FF}  |h_\epsilon(y)- h_0(y)|d\mu_{t_0,t_1,
  \dots,t_{l-1} }(y)\\
\leq 2\sup_\epsilon\sup_y|h_n(y)|\eta^{'} = 
2\sup_\epsilon\sup_y|\expect[g(\gamma_\epsilon(y))]|\eta^{'}.
\end{multline*} 
Thus $\transl{X}^\epsilon(t_l)\rightarrow \transl{X}^{0}(t_l)$ 
in distribution as $\epsilon\rightarrow0$.

We now prove by induction that 
$(\transl{X}_{t_0}^{\epsilon}, \transl{X}_{t_1}^{\epsilon}, 
                   \dots, \transl{X}_{t_{k}}^{\epsilon})$
converges in distribution to 
$(\transl{X}_{t_0}^{0}, \transl{X}_{t_1}^{0}, \dots,
\transl{X}_{t_{k}}^{0})$  
as $\epsilon\rightarrow0$.
By construction, we have 
$\transl{X}_{t_0}^{\epsilon}\rightarrow \transl{X}_{t_0}^{0}$ in distribution. 
Assume that for $0\leq l\leq k-1$, $(\transl{X}_{t_0}^{\epsilon},
\transl{X}_{t_1}^{\epsilon}, \dots, \transl{X}_{t_{l}}^{\epsilon})$ 
 converges in distribution in $\h_{2}^{l+1}$.
Let us define 
$
\alpha_\epsilon : \h_{2}^{l+1}\rightarrow  \h_{2}^{l+2}$
by 
$$\alpha_\epsilon(y_0, y_1, \dots, y_{l}) 
= \Bigl(y_0, y_1, \dots, y_{l}, S(t_{l+1}- t_l)y_{l} +
 \int_{t_{l}}^{t_{l+1}}S(t_{l+1}-s)F_\epsilon(s,y_{l})ds\Bigl)$$
and 
$
 \beta_\epsilon : \h_{2}^{l+1}\rightarrow  \ellp{1}(\esprob,\h_{2}^{l+2})
$
by 
$$\beta_\epsilon(y_0, y_1, \dots, y_{l}) 
= \Bigl(0,\dots, 0, \int_{t_{l}}^{t_{l+1}}S(t_{l+1}-s)G_\epsilon(s,y_l)dW(s)\Bigr)$$
so that 
$$\law{(\alpha_n +  \beta_n)
(\transl{X}_{t_0}^{n}, \transl{X}_{t_1}^{n}, \dots,\transl{X}_{t_{l}}^{n})} 
= 
\law{(\transl{X}_{t_0}^{n}, \transl{X}_{t_1}^{n}, \dots, 
 \transl{X}_{t_{l}}^{n},\transl{X}_{t_{l+1}}^{n} )}.$$
We denote
$u_\epsilon = 
(\transl{X}_{t_0}^{\epsilon}, \transl{X}_{t_1}^{\epsilon}, 
    \dots,\transl{X}_{t_{l}}^{\epsilon})$ 
and $\mu_\epsilon = \law{u_\epsilon}$.
Let $g\in \BL(\h_{2}^{l+2})$, and
$$h_\epsilon(y) = 
\expect g\bigl(\alpha_\epsilon(y) +  \beta_\epsilon(y)\bigl),\ 
y\in \h_{2}^{l+1}.$$
Proposition \ref{prop: Randomvar} yields
$$\expect g(
\transl{X}_{t_0}^{\epsilon}, \transl{X}_{t_1}^{\epsilon},
          \dots,\transl{X}_{t_{l+1}}^{\epsilon} ) 
= \expect h_\epsilon(u_\epsilon) = \int_{\h_{2}^{l+1}}h_\epsilon(y)d\mu_\epsilon(y).$$
It follows that
\begin{multline*}
|\expect g(\transl{X}_{t_0}^{\epsilon}, \transl{X}_{t_1}^{\epsilon}, 
      \dots,\transl{X}_{t_{l+1}}^{\epsilon} )- \expect g(\transl{X}_{t_0}^{0}, 
      \transl{X}_{t_1}^{0}, \dots,\transl{X}_{t_{l+1}}^{0} )|\\
\begin{aligned}
&\leq \int_{\h_{2}^{l+1}}|h_\epsilon(y)- h_0(y)|d\mu_\epsilon(y) 
+ |\int_{\h_{2}^{l+1}}h_0(y)d\mu_\epsilon(y) 
- \int_{\h_{2}^{l+1}}h_0(y)d\mu_0(y)|\\
&\leq J_1(\epsilon) + J_2(\epsilon). 
\end{aligned}
\end{multline*}
As in the above reasoning, we can prove that 
$$h_\epsilon\in \BL(\h_{2}^{l+1}),\, \,   \sup_\epsilon\|h_\epsilon\|_{\bl}<\infty$$ 
and 
$$\alpha_\epsilon(y) +  \beta_\epsilon(y) \rightarrow \alpha_0(y) +
\beta_0(y)$$ 
in distribution, for every  $y \in \h_{2}^{l+1}$
thus $h_\epsilon(y)\rightarrow h_0(y)$ as $\epsilon\rightarrow0$, 
for any $y \in \h_{2}^{l+1}$.
On the other hand, by the induction hypothesis we have 
$\mu_\epsilon\rightarrow\mu_0$  and since $h_0\in \BL(\h_{2}^{l+1})$
we have $ J_2(\epsilon)\rightarrow 0$ as $\epsilon\rightarrow0$. 
The convergence of $\mu_\epsilon$ implies that  $\{ \mu_\epsilon\}$ is tight,
i.e.~for each $\eta^{'}>0$  
there exists a compact set $\FF \in \h_{2}^{l+1}$ such that 
\begin{equation}\label{eq:Tight}
 \sup_\epsilon\mu_\epsilon(\h_{2}^{l+1}\setminus \FF) < \eta^{'}
\end{equation}
Since for every $y \in \h_{2}^{l+1}$, $h_\epsilon(y)\rightarrow h_0(y)$, 
from (\ref{eq:Tight}) and the compactness
 of $\FF$ the function $h_\epsilon$ converges to $h_0$ uniformly on $\FF$, hence 
 $$\lim_{\epsilon\rightarrow 0}\int_{\FF}  |h_\epsilon(y)-
 h_0(y)|d\mu_\epsilon(y) = 0$$ 
and
$$
\int_{\h_{2}^{l}\setminus\FF}  |h_\epsilon(y)-
h_0(y)|d\mu_\epsilon(y)\leq2\sup_\epsilon\sup_y|h_\epsilon(y)|\eta
$$
 So $ J_1(\epsilon)\rightarrow 0$ as $\epsilon\rightarrow0$, 
consequently 
$(\transl{X}_{t_0}^{\epsilon}, \transl{X}_{t_1}^{\epsilon},
\dots,\transl{X}_{t_{l+1}}^{\epsilon} )$ 
converges in distribution in $\h_{2}^{l+2}$.

Now, we have
\begin{multline*}
 d_{\bl}\Bigl(\law {Y^\epsilon(t_0),\dots,Y^\epsilon(t_k)}, 
             \law {Y^0(t_0),\dots,Y^0(t_k)}\Bigl)\\
\begin{aligned}
\leq & d_{\bl}\Bigl(
          \law{Y^\epsilon(t_0),\dots,Y^\epsilon(t_k)}, 
          \law{\transl{X}^\epsilon(t_0),\dots,\transl{X}^\epsilon(t_k)} \Bigr)\\
 & + d_{\bl}\Bigl(\law {\transl{X}^\epsilon(t_0),\dots,\transl{X}^\epsilon(t_k)}, 
                \law {\transl{X}^0(t_0),\dots,\transl{X}^0(t_k)}\Bigr)\\
 & + d_{\bl}\Bigl(\law {\transl{X}^0(t_0),\dots,\transl{X}^0(t_k)}, 
               \law {Y^0(t_0),\dots,Y^0(t_k)} \Bigr)\\
\leq &  \expect\Bigl(\lVert Y^\epsilon(t_0)-\transl{X}^\epsilon(t_0)\rVert+\dots
        + \lVert Y^\epsilon(t_k)-\transl{X}^\epsilon(t_k)\rVert\Bigr)\\
 &+d_{\bl}\Bigl(\law {\transl{X}^\epsilon(t_0),\dots,\transl{X}^\epsilon(t_k)}, 
                \law {\transl{X}^0(t_0),\dots,\transl{X}^0(t_k)}\Bigr)\\ 
& +\expect\Bigl(\lVert Y^0(t_0)-\transl{X}^0(t_0)\rVert+\dots
        + \lVert Y^0(t_k)-\transl{X}^0(t_k)\rVert\Bigr).
\end{aligned}
\end{multline*}
By 
\eqref{eq:partition_uniforme}, 
the first and third terms can be made
arbitrarily small. 
Thus  
$(Y^\epsilon(t_0),\dots,Y^\epsilon(t_k))$ converges in distribution to 
$(Y^0(t_0),\dots,Y^0(t_k))$. 
Now, for any finite sequence $\tau_1,\dots,\tau_m$  in $[\Tau,T]$, 
we can refine if necessary the partition $(t_1,\dots,t_k)$ such as to
include the points $\tau_1,\dots,\tau_m$.  
This proves that 
the 
finite dimensional distributions of 
$(Y^\epsilon)$
converge to the corresponding 
finite dimensional distributions of 
$Y^0.$

To show that $(Y^\epsilon)$ converges in distribution to $Y^0$ 
in $\C([\Tau,T], \h_2)$, we only need to prove that $(Y^\epsilon)$
is tight in $\C([\Tau,T], \h_2)$. 

By Condition (\ref{cond:pp}')-(a) 
and the equicontinuity of $(F_\epsilon)$, 
the sequence $\Bigl(\int_{\Tau}^{t}S(t-s)F_\epsilon(s,X(s))\,ds\Bigr)$ 
is tight. 
Tightness of $\Bigl(\int_{\Tau}^{t}S(t-s)G_\epsilon(s,X(s))\,dW(s)\Bigr)$
follows from  
\cite[Lemma 3.2]{jakubowski-kamenski-prf05stochinc}
applied to the multifunction
$$\GG(t,x)=\cobar\accol{G_\epsilon(t,x)\tq \epsilon>0},$$
(where $\cobar$ denotes the closed convex hull),
which is compact-valued thanks to Condition (\ref{cond:pp}')-(a). 
This result is given in \cite{jakubowski-kamenski-prf05stochinc} for
$p$-integrable stochastic processes with $p>2$, 
in order to use the stochastic convolution inequality. 
But the reasoning
remains unchanged for $p=2$ for a contractions semigroup, 
as is the case here,
because then the convolution equality still holds true, see 
\cite[Theorem 6.10]{dapratozabczyk92book}.



\medskip
\noindent{\em \textbf{Second step}.} 
We assume now that 
\begin{equation}\label{eq:XLp_bounded}
\sup_{t\in[\Tau,T]}\expect\norm{X(t)}^p<+\infty
\end{equation}
for some $p>2$.
Let us prove that  
 $\lim_{\epsilon\rightarrow0}\wass\CCO{Y^\varepsilon,Y^0}=0$ 
for any $T>0$.   
We only need to prove that 
$\CCO{\norm{Y^\epsilon}^2_{\C([\Tau,T],\h_2)}}$ is uniformly integrable, 
which is a consequence of \eqref{eq:XLp_bounded}:
we have, for any $\epsilon>0$,
\begin{align*}
\expect \CCO{\sup_{t\in[\Tau,T]}\lVert Y^\epsilon(t) \rVert^p}
\leq&  2^{p-1}\expect \CCO{\sup_{t\in[\Tau,T]}\lVert \int^t_{\Tau}S(t-s)
      F_\epsilon(s, X(s))\,ds\rVert^p}\\
    &+2^{p-1}\expect\CCO{\sup_{t\in[\Tau,T]} \lVert \int^t_{\Tau}S(t-s)
      G_\epsilon(s, X(s))\,dW(s)\rVert^p}\\
\leq& 2^{p-1}\expect \CCO{\sup_{t\in[\Tau,T]}\Bigr(\int^t_{\Tau}e^{-\delta(t-s)}
     \lVert  F_\epsilon(s, X(s))\rVert\,ds\Bigl)^p}\\
    &+2^{p-1}\convol\expect
               \int^T_{\Tau}\norm{G_\epsilon(s, X(s))}_{L(\h_1,
                 \h_2)}^{p}\,ds
\intertext{(where $\convol$ is given by the convolution inequality 
\cite[Proposition 7.3]{dapratozabczyk92book})}
\leq& 2^{p-1}\expect \CCO{\Bigl(\int^T_{\Tau}e^{-\delta(t-s)}
     \lVert  F_\epsilon(s, X(s))\rVert\,ds\Bigl)^p}\\
    &+2^{p-1}\convol\expect
               \int^T_{\Tau}\norm{G_\epsilon(s, X(s))}_{L(\h_1,
                 \h_2)}^{p}\,ds\\
\leq& 2^{p-1}(T-\Tau)^{p-1}\int^T_{\Tau}
     \lVert  F_\epsilon(s, X(s))\rVert^p\,ds\\
    &+2^{p-1}\convol\expect
               \int^T_{\Tau}\norm{G_\epsilon(s, X(s))}_{L(\h_1,
                 \h_2)}^{p}\,ds\\
\leq&2^{p-1}K^p\CCO{(T-\Tau)^{p-1}+\convol} 
                 \int^T_{\Tau} (1+\norm{X(s)})^p\,ds.  
\end{align*}
This estimation is independent of $\epsilon$, thus 
$(Y^\epsilon)_{\epsilon>0}$ is bounded in $\ellp{p}(\C([\Tau,T], \h_2))$. 
Thus 
$\CCO{\sup_{t\in[\Tau,T]}\lVert Y^\epsilon(t) \rVert^2}$ 
is uniformly integrable, 
which entails that $(Y^\epsilon)$ 
converges to $Y^0$ for $\wass$. 

It is then straightforward to deduce that, 
if $(X(t))_{t\in\R}$ is bounded in $\ellp{p}$, 
\begin{equation*}
\lim_{\epsilon\rightarrow0+}\wass(L_\epsilon(X),L_0(X))=0
\end{equation*}
for any interval $[\Tau,T]$
of $\R$.

\medskip
\noindent{\em \textbf{Third step}.} 
Now, let us show that  $L_\varepsilon(X^\varepsilon)$ converges to $L_0(X^0)$  
 in $\ellp{2}$-Wasserstein distance on compact intervals
 of $\R$  on the space $\C(\R, \h_2)$, which means that
$\wass(X^\varepsilon,X^0)\rightarrow0$ as 
 $\varepsilon\rightarrow 0+$ for any compact interval $[\Tau,T]$.

We have shown in the proof of Theorem \ref{theo:main} that
$L_\epsilon$ is $\theta$-Lipschitz for the norm of 
$\CUB\bigl(\R,\ellp{2}(\prob, \h_2)\bigl)$. 
This is not sufficient for our present purpose, 
but the hypothesis $\theta'<1$ allows a more precise
calculation: For $X,Y\in\CUB\bigl(\R,\ellp{2}(\prob, \h_2)\bigl)$, we have
\begin{multline*}
 \expect\CCO{\sup_{t\in\R}\| (L X)(t) -  (L Y)(t)\|_{\h_2}^2}\\
\begin{aligned}
\leq& 2 \expect\CCO{\sup_{t\in\R}\Bigl(\int_{-\infty}^{t}e^{-\delta(t-s)}
    \|F(s, X(s))- F(s,Y(s))\|_{\h_2}ds\Bigl)^2}\\
  & + 2\expect\CCO{\sup_{t\in\R}
\Bigl(\| \int_{-\infty}^{t}S(t-s)[G(s, X(s)) - G(s,Y(s))]dW(s)\|_{\h_2}\Bigl)^2}\\
  = &  J_1 + J_2.
\end{aligned}
\end{multline*}
We have
\begin{align*}
J_1   
\leq& 2 \expect\CCO{\sup_{t\in\R}\Bigl(\int_{-\infty}^{t}e^{-\delta(t-s)}
   \sup_{\sigma\in\R} 
       \|F(\sigma, X(\sigma))- F(\sigma,Y(\sigma))\|_{\h_2}ds\Bigl)^2}
\\
=&2\expect\CCO{\sup_{\sigma\in\R} 
       \|F(\sigma, X(\sigma))- F(\sigma,Y(\sigma))\|_{\h_2}^2}
\sup_{t\in\R}\Bigl(\int_{-\infty}^{t}e^{-\delta(t-s)}ds\Bigl)^2\\
\leq& \frac{2}{\delta^2}
\expect\CCO{\sup_{t\in\R}
\|F(t, X(t))-F(t,Y(t))\|_{\h_2}^2ds}\\
\leq& \frac{2K^2}{\delta^2}
          \sup_{t\in\R}\expect\|X(t))- Y(t))\|_{\h_2}^2.
\end{align*}
By \eqref{cond:croissance} and \eqref{cond:lipschitz}, the process 
$\int_\infty^.S(.-s)\CCO{G(s,X(s))-G(s,Y(s))}\,dW(s)$ is a square
integrable martingale. Using Doob's inequality and Itô's isometry
identity, we get
\begin{align*}
J_2 =&2\sup_{T\in\R}
\expect\CCO{\sup_{t\leq T}
\Bigl(\| \int_{-\infty}^{t}S(t-s)[G(s, X(s)) - G(s,Y(s))]\,dW(s)\|_{\h_2}\Bigl)^2}\\
\leq& 8\sup_{T\in\R}\expect
\Bigl(\| \int_{-\infty}^{T}S(T-s)[G(s, X(s)) - G(s,Y(s))]\,dW(s)\|_{\h_2}\Bigl)^2\\
\leq& 8\sup_{T\in\R}\trace Q 
   \int_{-\infty}^{T}e^{-2\delta(T-s)}
         \expect\|G(s, X(s)) - G(s,Y(s))\|^{2}_{L(\h_1, \h_2)}ds\\
\leq& 8\sup_{T\in\R}\trace QK^2
    \int_{-\infty}^{T}e^{-2\delta(T-s)}
      \expect\|X(s) - Y(s)\|^{2}_{\h_2}ds\\
\leq &8 \sup_{T\in\R}K^2\trace Q\bigl(\int_{-\infty}^{T}e^{-2\delta(t-s)}ds\bigl)\,
                   \sup_{\sigma\in\R} \expect\|X(\sigma) - Y(\sigma)\|^{2}_{\h_2}\\
\leq&\dfrac{4K^2\trace Q}{\delta}\sup_{t\in\R} \expect\|X(t) - Y(t)\|^{2}_{\h_2}.
\end{align*}
Thus 
\begin{multline*}
\expect\CCO{\sup_{t\in\R} \|(LX)(t) -  (LY)(t)\|_{\h_2}^2 }\leq J_1 + J_2 \\
           \leq {\frac{2K^2}{\delta}\CCO{\frac{1}{\delta}+2\trace{Q}}} 
                \expect\CCO{\sup_{t\in\R} \|X(t) - Y(t)\|^{2}_{\h_2}}.
\end{multline*}
with 
$$\frac{2K^2}{\delta}\CCO{\frac{1}{\delta}+2\trace{Q}}<
\dfrac{4K^2}{\delta}\CCO{\dfrac{1}{\delta} + \trace Q}=\theta'.$$
We deduce
\begin{align*}
\wass(X^\epsilon,X^0)
= &\wass(L_\epsilon(X^\epsilon),L_0(X^0))\\
\leq & \wass(L_\epsilon(X^\epsilon),L_\epsilon(X^0))
      +\wass(L_\epsilon(X^0),L_0(X^0))\\
\leq &\theta'\wass(X^\epsilon,X^0)
      +\wass(L_\epsilon(X^0),L_0(X^0)).
\end{align*}
As $\theta'<1$, this entails
\begin{equation}\label{eq:wass2zero}
\wass(X^\epsilon,X^0)\leq
\frac{1}{1-\theta'}\wass(L_\epsilon(X^0),L_0(X^0)). 
\end{equation}
But, as $\theta'<1$, there exists $p>2$ such that 
$(X^0(t))_{t\in\R}$ is bounded in $\ellp{p}$, 
see Remark \ref{rem:novikov-BDG} and 
Proposition \ref{cor:Lpboundedness}.
Thus the right hand side of \eqref{eq:wass2zero}
converges to $0$
as $\epsilon\rightarrow 0+$.

Finally, by \cite[Theorem 4.1]{DaPrato-Tudor95}, 
the mild solution  $(X^0)$  to (\ref{eq:station}) is 
 a stationary process. 
\finpr

\paragraph{Acknowledgements} 
We are grateful to Professor Fazia Bedouhene, 
from the university of Tizi Ouzou,
and Professor 
Zhenxin Liu, from Jilin University, who independently 
pointed out to us errors in earlier versions of the proof of 
Theorem \ref{theo:main}.

\def\cprime{$'$} \def\cprime{$'$} \def\cprime{$'$} \def\cprime{$'$}

\end{document}